\documentclass{amsart}
\font\bit=cmbxti10                       
\font\sbit=cmbxti10 scaled 700           

\def\bb{{\hbox{\bit b}}}

\newsymbol\diagdown 231F
\newsymbol\digamma 207A
\newsymbol\varnothing 203F

\let\void=\varnothing
\let\<=\langle
\let\>=\rangle
\let\veps=\varepsilon

\let\what=\widehat
\let\\=\backslash
\let\sse=\subseteq
\let\limply=\Longrightarrow

\def\newmatrix#1{\null\,\vcenter{
		\baselineskip=8pt\mathsurround=-0pt\ialign{
		\hfil ${##}$
		\hfil &&
		\hfil ${##}$
		\hfil \crcr
		\mathstrut \crcr
		\noalign{\kern-\baselineskip}#1 \crcr
		\mathstrut \crcr
		\noalign{\kern-\baselineskip} \crcr }}\!}

\def\CC{{\mathbb C\kern.5pt}}
\def\FF{{\mathbb F\kern.5pt}}
\def\RR{{\mathbb R\kern.5pt}}
\def\B{{\mathcal B}}
\def\F{{\mathcal F}}
\def\H{{\mathcal H}}
\def\J{{\mathcal J}}
\def\K{{\mathcal K}}
\def\Le{{\mathcal L}}

\def\N{{\mathcal N}}
\def\R{{\mathcal R}}

\def\X{{\mathcal X}}
\def\Y{{\mathcal Y}}
\def\V{{\kern-.5pt\mathcal V}}
\def\W{{\mathcal W}}
\def\Z{{\mathcal Z}}

\def\ssb{{\kern-.5pt\hbox{\sbit b}}}
\def\ssB{{\kern-.5pt\scriptscriptstyle\B}}

\def\noi{\noindent}
\def\tr{{\rm tr}}
\def\0{{\{0\}}}
\def\half{{{\frac{1}{2}}}}
\def\smallfrac#1#2{{\textstyle{\frac{#1}{#2}}}}
\def\span{{\kern.5pt{\rm span}\kern2pt}}

\def\x{{\times}}
\def\ve{{_\vee\kern-1pt}}
\def\we{{_\wedge\kern-1pt}}
\def\hotimes{{\kern1pt\what\otimes}}
\def\otimesv{{\kern1pt\otimes_\ve}}
\def\otimesw{{\kern1pt\otimes_{\kern-1pt\we}}}
\def\hotimesv{{\kern1pt\what\otimes_\ve}}
\def\hotimesw{{\kern1pt\what\otimes_{\kern-1pt_\we}}}

\def\emap{\hbox to25pt{\rightarrowfill}}

\def\nmap{\Big\uparrow}
\def\diagdownBOX{\hbox{$\diagdown$}}
\def\searrowBOX{\hbox{\hglue6.5pt$\searrow$}}
\def\semap{\vbox{\offinterlineskip
\diagdownBOX\vglue-1pt\searrowBOX\vglue-6pt}}

\def\theorem #1{\vskip6pt\noi{\bf{Theorem} #1$.$}}
\def\lemma #1{\vskip6pt\noi{\bf{Lemma} #1$.$}}
\def\corollary #1{\vskip6pt\noi{\bf{Corollary} #1$.$}}
\def\proposition #1{\vskip6pt\noi{\bf{Proposition} #1$.$}}
\def\claim #1{\vskip6pt\noi{\it{\kern-.5pt Claim} #1$.$}}
\def\remark #1{\vskip6pt\noi{\bf{Remark} #1$.$}}

\begin{document}

\vglue-55pt\noi
\hfill{\it Concrete Operators}\/,
{\bf 9} (2022) 53--69

\vglue20pt
\title
[Trace-Class and Nuclear Operators]
{Trace-Class and Nuclear Operators}
\author{Carlos S. Kubrusly}
\address{Catholic University of Rio de Janeiro, Brazil}
\email{carlos@ele.puc-rio.br}
\renewcommand{\keywordsname}{Keywords}
\date{December 16, 2021; revised$:$ April 11, 2022}

\begin{abstract}
This paper explores the long journey from projective tensor products of a pair
of Banach spaces, passing through the definition of nuclear operators still on
the realm of projective tensor products, to the of notion of trace-class
operators on a Hilbert space, and shows how and why these concepts (nuclear
and trace-class operators, that is) agree in the end.
\end{abstract}

\maketitle

\vskip-20pt\noi
\section{Introduction}
This is an expository paper on trace-class and nuclear operators$.$ Its
purpose is to demonstrate that these classes of operators coincide on a
Hilbert space$.$ It will focus mainly on three points: (i) where these
notions came from, (ii) how they are intertwined, and (iii) when they
coincide$.$ Being a plain expository paper, this has no intention to survey
the subject, neither to offer an extensive bibliography on it$.$

\vskip6pt
To begin with, let us borrow the definitions of nuclear and trace-class
operators from Sections 4 and 5 (where these definitions will be properly
posed).
\begin{description}
\item{$\circ$}
{\sc Nuclear Operators}$.$
An operator $T$ on a Banach space $\X$ is {\it nuclear}\/ if there are
$\X^*\!$-valued and $\X$-valued sequences $\{f_k\}$ and $\{y_k\}$ such that
${\sum_k\|f_k\|\kern1pt\|y_k\|}<\infty$ and $Tx={\sum_k f_k(x)\kern1pt y_k}$
for every ${x\in\X}$.
\vskip6pt
\item{$\circ$}
{\sc Trace-Class Operators}$.$
An operator $T$ on a Hilbert space $\X$ is {\it trace-class}\/ if
${\sum_\gamma\!\big\<|T|e_\gamma;e_\gamma\big\><\infty}$ for an arbitrary
orthonormal basis $\{e_\gamma\}$ for $\X$ and the series value does not
depend on $\{e_\gamma\}$.
\end{description}

\vskip6pt
Perhaps the first lines in Robert Bartle's review of \cite{Ha} might be a
suitable start$:$ ``Grothendieck \cite{Gro} showed that a Banach space $\X$
has the approximation property if and only if, for every nuclear operator
${T\!:\X\!\to\X}$ (i.e., operator having the form
$T=$ ${\sum_k\<f_k\,;\cdot\>y_k}$) with ${f_k\in\X^*}\!$, $\,{y_k\in\X}$, and
${\sum_k\|f_k\|\,\|y_n\|<\infty}$), the number $\tr(T)=$ ${\sum_k\<f_k,y_k\>}$
is well-defined (i.e., is independent of the choice of $\{f_k\}$ and $\{y_k\}$
in the representation $T={\sum_k\<f_k\,;\cdot\>y_k}$) and can be used to
define a trace function.''

\vskip6pt
Bartle's concise description nicely summarizes the apparently long way to be
covered from Grothendieck's projective tensor products (where nuclear
operators originate form), to trace-class operators, and finally concluding
that these classes coincide$.$ The familiar notion of trace as the sum of
eigenvalues is a fundamental result known as Lidski\v{\i} Theorem \cite{Lid},
\cite[Theorem 8.4]{GK} which still remains an active research topic (e.g.,
\cite{Pow,JS,FJ,Did,Rei,Pie}).

\vskip6pt
These concepts (nuclear and trace-class) are linked together since their early
days$.$ Schatten in his celebrated 1950 monograph \cite{Sch1} (which actually
is an offspring of his 1942 thesis) describes nuclear Hilbert-space operators
as being precisely the trace-class$:$ ``The trace-class may be also
interpreted as ${\X\hotimesw\X^*}$'' \cite[Theorem 5.12]{Sch1} --- the
completion of the tensor product of a Hilbert with its dual, with respect to
the greatest reasonable crossnorm (i.e., the projective norm).

\vskip6pt
The paper is organized as follows$.$ Section 2 summarizes some common notation
and terminology$.$ Section 3 poses the necessary results on crossnormed
tensor products of a pair of Banach spaces, since this is the proper setup
where nuclear transformations come from$.$ Nuclear transformations are defined
in Section 4 (Theorem 4.1) as the range of a linear contraction
of the projective tensor product ${\X^*\hotimesw\Y}$ into the Banach space
${\B[\X,\Y]}$ of all bounded linear transformations from $\X$ to $\Y$, which
yields the characterization of nuclear operators mentioned above$.$ Section 5
deals exclusively with trace-class operators on a Hilbert space, and gives a
thorough view of basic properties of these operators$.$ Section 6 shows in
Theorem 6.1 that nuclear and trace-class operators in fact reduce to the
same thing.

\vskip6pt
All terms and notation above will be defined here in due course.

\section{Notation and Terminology}

Throughout the paper all linear spaces are over the same field $\FF$ (and the
field $\FF$ in this context means either $\RR$ or $\CC).$ If ${\X,\Y,\Z}$ are
linear spaces, then let ${\Le[\X,\Z]}$ and ${b[\X\x\Y,\Z]}$ denote the linear
spaces of all linear transformation of $\X$ into $\Z$, and of all bilinear
maps of the Cartesian product ${\X\x\Y}$ into $\Z.$ If ${\X,\Y,\Z}$ are normed
spaces, then let ${\B[\X,\Z]}$ denote the normed space of all bounded (i.e.,
continuous) linear transformations of $\X$ into $\Z$ equipped with its standard
induced uniform norm, and let ${\bb[\X\x\Y,\Z]}$ denote the normed spaces of
all bounded (i.e., continuous) bilinear maps of ${\X\x\Y}$ into $\Z$ equipped
with its usual norm, which are both Banach spaces whenever $\Z$ is (see e.g.,
\cite{Kub3}).

\vskip6pt
A subspace of a normed space is a closed linear manifold of it$.$ If $\X$ is
a normed space, then ${\X^*\!=\B[\X,\FF]}$ stand for its dual; if $M$ is a
subset of an inner product space, then $M^\perp$ stands for the orthogonal
complement of $M.$ Range and kernel of a bounded linear transformation
${T\kern-1pt\in\B[\X,\Y]}$ between normed spaces $\X$ and $\Y$ will be denoted
by $\R(T)$ --- a linear manifold of $\Y$ --- and $\N(T)$ --- a subspace of
$\X$ --- respectively$.$ If two normed spaces $\X$ and $\Y$ are isometrically
isomorphic, and if ${y\in\Y}$ is the isometrically isomorphic image of
${x\in\X}$, then write ${\X\cong\Y}$ and ${x\cong y}$.

\vskip6pt
Let $\X$ and $\Y$ be arbitrary normed spaces$.$ By an operator we mean a
bounded linear transformation of a normed space into itself$.$ Set
$\B[\X]={\B[\X,\X]}$; the normed algebra of all operators on $\X.$ Let
${\B_0[\X,\Y]}$ and ${\B_\infty[\X,\Y]}$ stand for the normed spaces of all
bounded finite-rank (i.e., ${\dim\R(T)<\infty}$) and of all compact linear
transformations ${T\!:\X\to\Y}$, respectively$.$ Similarly, set
$\B_0[\X]={\B_0[\X,\X]}$ and $\B_\infty[\X]={\B_\infty[\X,\X]}$; the ideals of
the algebra $\B[\X]$ consisting of all bounded finite-rank and of all
compact operators, respectively, so that
${\B_0[\X]\sse\B_\infty[\X]\sse\B[\X]}.$ Moreover, let ${\B_N[\X,\Y]}$
(accordingly, ${\B_N[\X]=B_N[\X,\X]}$) stand for the normed spaces of all
nuclear transformations$.$ It is worth noticing that there are different, also
common, notations such as $\K$ for compact, $\F$ for finite-rank, and $\N$ for
nuclear operators --- at the end of Section 5 it will become clear the reason
for our choice in favor of the above sub-indexed-$\B$ notation (see also,
e.g., \cite[Sections 6.1 and 7.1]{Wei})$.$ A Banach space $\Y$ has the
{\it approximation property}\/ if ${\B_0[\X,\Y]}$ is dense in
${\B_\infty[\X,\Y]}$ for every normed space $\X$ --- every Banach space with a
Schauder basis has the approximation property, in particular, since the range
of a compact linear transformation is separable, every Hilbert space has the
approximation property.

\section{Preliminaries on Crossnormed Tensor Product Spaces}

The algebraic {\it tensor product}\/ of linear spaces $\X$ and $\Y$ is a
linear space ${\X\otimes\Y}$ for which there is a bilinear map
${\theta\!:\X\x\Y\to\X\otimes\Y}$ (called the {\it natural bilinear map}\/
associated with ${\X\otimes\Y}$) whose range spans ${\X\otimes\Y}$ with the
following additional property: for every bilinear map ${\phi\!:\X\x\Y\to\Z}$
into any linear space $\Z$ there exists a (unique) linear transformation
${\Phi\!:\X\otimes\Y\to\Z}$ for which the diagram
$$
\newmatrix{
\X\x\Y & \kern2pt\buildrel\phi\over\emap & \kern-1pt\Z                   \cr
       &                                 &                               \cr
       & \kern-3pt_\theta\kern-3pt\semap & \kern4pt\nmap\scriptstyle\Phi \cr
       &                                 & \phantom{;}                   \cr
       &                                 & \kern-2pt\X\otimes\Y          \cr}
$$
commutes$.$ Set ${x\otimes y=\theta(x,y)}$ for each ${(x,y)\in\X\x\Y}.$ These
are the {\it single tensors}\/$.$ An arbitrary element $\digamma$ in the
linear space ${\X\otimes\Y}$ is a finite sum ${\sum_ix_i\otimes y_i}$ of
single tensors, and the representation of $\digamma$ as a finite sum of single
tensors ${\digamma=\sum_ix_i\otimes y_i}$ is not unique$.$ (For an exposition
on algebraic tensor products see, e.g., \cite{Kub2}.)

\vskip6pt
If $\X$ and $\Y$ are Banach spaces and $\X^*$ and $\Y^*$ are their duals, then
let ${x\otimes y}$ and ${f\otimes g}$ be single tensors in the tensor product
spaces ${\X\otimes\Y}$ and ${\X^*\otimes\Y^*}.$ A norm ${\|\cdot\|_\alpha}$ on
${\X\otimes\Y}$ is a {\it reasonable crossnorm}\/ if, for every
${x\kern-1pt\in\kern-1pt\X}$, ${y\kern-1pt\in\kern-1pt\Y}$,
${f\kern-1pt\in\kern-1pt\X^*}\!\!$, ${g\kern-1pt\in\kern-1pt\Y^*}\!\!$,
\begin{description}
\item {$\kern-12pt$\rm(a)}
$\;\|x\otimes y\|_\alpha\le\|x\|\,\|y\|$,
\vskip4pt
\item {$\kern-12pt$\rm(b)}
$\;{f\otimes g}$ lies in $(\X\otimes\Y)^*\!$, $\,$and
$\|f\otimes g\|_{*\alpha}\le\|f\|\,\|g\|$
\qquad
(where ${\|\cdot\|_{*\alpha}}$ is the norm on the dual $(\X\otimes\Y)^*$ when
$(\X\otimes\Y)$ is equipped with the norm ${\|\cdot\|_\alpha}$),
\end{description}
so that ${\X^*\otimes\Y^*\sse(\X\otimes\Y)^*}\!.$ It can be verified that (i)
${\|x\otimes y\|_\alpha}={\|x\|\,\|y\|}$ whenever ${\|\cdot\|_\alpha}$ is a
reasonable crossnorm, and (ii) when restricted to ${\X^*\otimes\Y^*}$ the norm
${\|\cdot\|_{*\alpha}}$ on ${(\X\otimes\Y)^*}$ is a reasonable crossnorm (with
respect to ${(\X^*\otimes\Y^*)^*}).$ Two special norms on ${\X\otimes\Y}$ are
the so-called {\it injective}\/ ${\|\cdot\|_\ve}$ and {\it projective}\/
${\|\cdot\|_\we}$ norms,
$$
\|\digamma\|_{_\vee\kern-1pt}
=\sup_{\|f\|\le1,\,\|g\|\le1,\;f\in\X^*\!,\,g\in\Y*}
\Big|{\sum}_if(x_i)\,g(y_i)\Big|,
$$
$$
\|\digamma\|_{_\wedge\kern-1pt}
=\inf_{\{x_i\}_i,\,\{y_i\}_i,\;\digamma=\sum_ix_i\otimes y_i}
{\sum}_i\|x_i\|\,\|y_i\|,
$$
for every ${\digamma=\sum_ix_i\otimes y_i\in\X\otimes\Y}$, where the infimum
is taken over all (finite) representations of ${\digamma\in\X\otimes\Y}.$
It can also be shown that (iii) these are indeed norms on ${\X\otimes\Y}$,
that (iv) both ${\|\cdot\|_\ve}$ and ${\|\cdot\|_\we}$ are reasonable
crossnorms and, moreover, that (v) a norm ${\|\cdot\|_\alpha}$ on
${\X\otimes\Y}$ is a reasonable crossnorm if and only if
$$
\|\digamma\|_{_\vee\kern-1pt}
\le\|\digamma\|_\alpha
\le\|\digamma\|_{_\wedge\kern-1pt}
\quad\hbox{for every}\;\quad
\digamma\in\X\otimes\Y.
$$
Let ${\X\otimes_\alpha\!\Y=(\X\otimes\Y,\|\cdot\|_\alpha)}$ stand for the
tensor product space of a pair of \hbox{Banach} spaces equipped with a
reasonable crossnorm ${\|\cdot\|_\alpha}$, which is not necessarily
com\-plete$.$ Their completion (see, e.g., \cite[Section 4.7]{EOT}) is denoted
by ${\X\hotimes_\alpha\Y}$ (same notation ${\|\cdot\|_\alpha}$ for the
extended norm on ${\X\hotimes_\alpha\Y}$)$.$ In particular, ${\X\hotimesv\Y}$
and ${\X\hotimesw\Y}$ are referred to as the {\it injective and projective
tensor products}\/$.$ For the theory of the Banach space
${\X\hotimes_\alpha\Y}$ (in particular, ${\X\hotimesv\Y}$ and
${\X\hotimesw\Y}$) see, e.g., \cite{Jar}, \cite{DF}, \cite{Rya}, \cite{DFS}.

\vskip6pt
The next two fundamental results on the projective tensor product will be
need along the next two sections$.$ The first one is referred to as
Grothendieck Theorem (see, e.g., \cite[Proposition 1.1.4]{DFS} and
\cite[Proposition 2.8]{Rya})$.$ In fact, most results mentioned in this
section (and beyond) are Grothendieck's, and there are different theorems
named after Grothendieck (the one in Theorem 3.1 below is not the classi\-cal
Grothendieck Theorem as in, for instance, \cite[Chapter 4]{DFS} and
\cite{Pis2})$.$ The second result in Theorem 3.2 below is called the Universal
Mapping Principle (see, e,g., \cite[Theorem 1.1.8]{DFS} and
\cite[Theorem 2.9]{Rya})$.$ Proofs are included for sake of completeness.

\theorem{3.1} {\sc (Grothendieck)}$.$
{\it If\/ $\X$ and\/ $\Y$ are Banach spaces, then for every\/
${\digamma\in\X\hotimesw\Y}$ there exist $\X$-valued and $\Y$-valued sequences
$\{x_k\}$ and $\{y_k\}$, respectively, for which the real sequence
$\{\|x_k\|\,\|y_k\|\}$ is summable and
$$
\digamma={\sum}_kx_k\otimes y_k
$$
$($i.e., ${\digamma\!\in\X\hotimesw\Y}$ is representable in the form of
a countable sum\/ $\digamma\!={\sum_kx_k\kern-1pt\otimes y_k}$ in the
sense that it is either a finite or countably infinite repre\-sentation$).$
Moreover, the projective norm ${\|\cdot\|_\we}$ on ${\X\hotimesw\Y}$ is given
by
$$
\|\digamma\|_\we=\inf{\sum}_k\|x_k\|\,\|y_k\|,
$$
where the infimum is taken over all representations\/
${\sum_kx_k\kern-1pt\otimes y_k}$ of}\/ ${\digamma\kern-1pt\in\X\hotimesw\Y}$.

\proof
Let ${\X\hotimesw\Y}$ be the completion of ${\X\otimesw\Y}.$ Identify
${\X\otimesw\Y}$ with isometrically isomorphic images of it, and so regard
${\X\otimesw\Y}$ as being dense in ${\X\hotimesw\Y}.$ Thus take
${\digamma\!\in\kern-1pt\X\hotimesw\Y\,\\\X\otimesw\Y}$ (otherwise the
resulting finite sum is trivially obtained) so that $\digamma$ is arbitrarily
close to elements in ${\X\otimesw\Y}.$ Take an arbitrary ${\veps>0}.$ For
each posi\-tive integer $k$ take
$\digamma_k\!={\sum_{i=1}^{n_k}x_i\otimes y_i}$ in ${\X\otimesw\Y}$ such that
$$
\|\digamma-\digamma_k\|_\we<\,\smallfrac{\veps}{2}\smallfrac{1}{2^k}.
$$
In particular,
$\big|\kern1pt\|\digamma_1\|_\we-\|\digamma\|_\we\big|
\le\|\digamma_1-\digamma\|_\we\le\frac{\veps}{2}\frac{1}{2}$
and hence $\|\digamma_1\|_\we\le\|\digamma\|_\we+\frac{\veps}{2}\frac{1}{2}.$
Take a representation ${\sum_{i=1}^{n_1}x_i\otimes y_i}$ for $\digamma_1$ for
which ${\sum_{i=1}^{n_1}\|x_i\|\,\|y_i\|}$ is close enough to
${\|\digamma_1\|_\we}
=\inf_{\{x_i\}_{i=1}^{n_1},\,\{y_i\}_{i=1}^{n_1}}
{\sum_{i=1}^{n_1}\|x_i\|\,\|y_i\|}$,
say
$$
{\sum}_{i=1}^{n_1}\|x_i\|\,\|y_i\|
\le\|\digamma_1\|_\we\!+\smallfrac{3\veps}{4}
\le\|\digamma\|_\we\!+\veps.
$$
Computing the norm of ${\digamma_{k+1}-\digamma_k\in\X\otimesw\Y}$ we get
$$
\|\digamma_{k+1}-\digamma_k\|_\we
\le\|\digamma_{k+1}-\digamma\|_\we\!+\|\digamma-\digamma_k\|_\we
<\smallfrac{\veps}{2}\smallfrac{1}{2^{k+1}}
+\smallfrac{\veps}{2}\smallfrac{1}{2^k}
=\smallfrac{\veps}{2}\smallfrac{1}{2^k}\smallfrac{3}{2}
<\veps\smallfrac{1}{2^k},
$$
so we can take a representation
${\sum_{i=n_k+1}^{n_{k+1}}x_i\otimes y_i}$ for
${\digamma_{k+1}-\digamma_k\in\X\otimesv\Y}$ for which
${\sum_{i=n_k+1}^{n_{k+1}}\|x_i\|\,\|y_i\|}$ is close enough to
${\|\digamma_{k+1}-\digamma_k}\|_\we$, say
$$
{\sum}_{i=n_k+1}^{n_{k+1}}\|x_i\|\,\|y_i\|
\le\|\digamma_{k+1}-\digamma_k\|_\we\!+\smallfrac{\veps}{4}\smallfrac{1}{2^k}
<\smallfrac{\veps}{2}\smallfrac{1}{2^k}\smallfrac{3}{2}
+\smallfrac{\veps}{2}\smallfrac{1}{2^k}\smallfrac{1}{2}
=\veps\smallfrac{1}{2^k}.
$$
Since ${\|\cdot\|_\we}$ is a reasonable crossnorm we get for each ${k>1}$
\begin{eqnarray*}
{\sum}_{i=1}^{n_{k+1}}\|x_i\otimes y_i\|_\we
&=&
{\sum}_{i=1}^{n_1}\|x_i\|\,\|y_i\|
+{\sum}_{j=1}^k{\sum}_{i=n_j+1}^{n_{j+1}}\|x_i\|\,\|y_i\|                  \\
&<&
\|\digamma\|_\we\!+\veps+\veps{\sum}_{j=1}^k\smallfrac{1}{2^j}
<\|\digamma\|_\we\!+2\veps.
\end{eqnarray*}
Thus the sequence $\{x_i\otimes y_i\}$ is absolutely summable with
${\sum_i\kern-1pt\|x_i\otimes y_i\|_\we}\!<{\|\digamma\|_\we\!+2\veps}.$
Therefore, since ${\X\hotimesw\Y}$ is a Banach space, the sequence is
summable$.$ This means the sequence of partial sums
$\{{\sum_{i=1}^{n_k}x_i\otimes y_i}\}$ converges in ${\X\hotimesw\Y}.$ Hence
the limit ${\digamma\!\in\X\hotimesw\Y}$ of
$\digamma_k\!
={\sum_{i=1}^{n_k}x_i\otimes y_i\in\X\otimesw\Y}$ is represented as
$$
\digamma={\sum}_ix_i\otimes y_i,
$$
\vskip-6pt\noi
and
$$
\|\digamma\|_\we
\le{\sum}_i\|x_i\otimes y_i\|_\we
={\sum}_i\|x_i\|\,\|y_i\|
\le\|\digamma\|_\we+2\veps,
$$
so that the projective norm $\|\digamma\|_\we$ of $\digamma$ is the infimum
of ${\sum_i\|x_i\|\,\|y_i\|}$ over all representations of $\digamma$
of the form ${\sum_ix_i\otimes y_i}.$ \qed

\theorem{3.2} {\sc Universal Mapping Principle}$.$
{\it If\/ ${(\X,\Y,\Z)}$ is an arbitrary triple of Banach spaces, then
$$
\bb[\X\x\Y,\Z]\cong\B[\X\hotimesw\Y,\Z]
$$
$($i.e., the Banach spaces\/ ${\bb[\X\x\Y,\Z]}$ and\/
${\B[\X\hotimesw\Y,\Z]}$ are isometrically isomorphic}\/).

\proof
Take the natural bilinear map ${\theta\in b[\X\x\Y,\X\otimes\Y]}$ associated
with the tensor product space ${\X\otimes\Y}.$ It is readily verified that
the composition with $\theta$, namely
$$
C_\theta(\Phi)=\Phi\circ\theta\,\in\,b[\X\x\Y,\Z]
\quad\;\hbox{for every}\,\quad
\Phi\in\Le[\X\otimes\Y,\Z],
$$
defines a linear-space isomorphism
${C_\theta\!:\Le[\X\otimes\Y,\Z]\to b[\X\x\Y,\Z]}.$ Equip the linear space
${\X\otimes\Y}$ with the projective norm to get the normed space
${\X\otimesw\Y}.$ Since $\Z$ is a Banach space, ${\B[\X\otimesw\Y,\Z]}$
is a Banach space which is a linear manifold of the linear space
${\Le[\X\otimesw\Y,\Z]}$. Let $\J$ be the restriction to
${\B[\X\otimesw\Y,\Z]}$ of the linear-space isomorphism $C_\theta$ on
${\Le[\X\otimesw\Y,\Z]}$, which remains linear and injective,
$$
\J=C_\theta|_{\B[\X\otimesw\Y,\Z]}
\!:\B[\X\otimesw\Y,\Z]\to\R(\J)\sse b[\X\x\Y,\Z].
$$
Next we show that (a) ${\R(\J)=\bb[\X\x\Y,\Z]}$ and (b) $\J$ is an isometry
onto $\R(\J)$. (Let ${\|\cdot\|_{_{_\ssB}}}$ and ${\|\cdot\|_{_\ssb}}$ stand
for the norms in ${\B[\X\otimesw\Y,\Z]}$ and ${\bb[\X\x\Y,\Z]}$.)

\vskip6pt\noi
(a$_1$)
If ${\Phi\in\B[\X\otimesw\Y,\Z]}$, then
${\J(\Phi)=C_\theta(\Phi)=\Phi\circ\theta\in b[\X\x\Y,\Z]}$ is such that
$$
\phi(x,y)=\J(\Phi)(x,y)=(\Phi\circ\theta)(x,y)=\Phi(x\otimes y)
\quad\;\hbox{for every}\;\quad
(x,y)\in\X\x\Y,
$$
and hence
$$
\|\J(\Phi)(x,y)\|_{_\Z}
=\|\Phi(x\otimes y)\|_{_\Z}
\le\|\Phi\|_{_{_\ssB}}\|x\otimes y\|_{_\we}
=\|\Phi\|_{_{_\ssB}}\|x\|_{_\X}\|y\|_{_\Y}.                    \eqno{(*)}
$$
Then the bilinear map $\J(\Phi)$ is bounded. Thus
${\J(\Phi)\in\bb[\X\x\Y,\Z]}$ for every $\Phi$ in ${\B[\X\otimesw\Y,\Z]}$,
and so ${\R(\J)\sse\bb[\X\x\Y,\Z]}$.

\vskip6pt\noi
(a$_2$)
Conversely, if ${\phi\in\bb[\X\x\Y,\Z]}$, then there is a unique
${\Phi\in\B[\X\otimesw\Y,\Z]}$ such that
${C_\theta(\Phi)=\Phi\circ\theta}=\phi$. So
$$
\|\Phi(x\otimes y)\|_{_\Z}
=\|\Phi\circ\theta(x,y)\|_{_\Z}
=\|\phi(x,y)\|_{_\Z}
\le\|\phi\|_{_{\ssb}}\|x\|_{_\X}\|y\|_{_\Y}
$$
for every ${(x,y)\in\X\x\Y}$. Hence for an arbitrary
${\digamma=\sum_ix_i\otimes y_i\in\X\otimesw\Y}$,
\begin{eqnarray*}
\|\Phi(\digamma)\|_{_\Z}
&=&
\Big\|\Phi{\sum}_ix_i\otimes y_i\Big\|_{\scriptscriptstyle\Z}
=\Big\|{\sum}_i\Phi(x_i\otimes y_i)\Big\|_{\scriptscriptstyle\Z}           \\
&\le&
{\sum}_i\|\phi\|_{_\ssb}\|x_i\|_{_\X}\|y_i\|_{_\Y}
=\|\phi\|_{_\ssb}{\sum}_i\|x_i\|_{_\X}\|y_i\|_{_\Y}.
\end{eqnarray*}
Since this holds for every (finite) representation of
$\digamma\!=\sum_ix_i\otimes y_i$, and since
$\|\digamma\|_\we\!=\inf{\sum}_i\,\|x_i\|_{_\X}\|y_i\|_{_\Y}$ over all
representations, then for every ${\digamma\in\X\otimesw\Y}$
$$
\|\Phi(\digamma)\|_{_\Z}\le\|\phi\|_{_\ssb}\|\digamma\|_\we,    \eqno{(**)}
$$
and the linear $\Phi$ is bounded: ${\Phi\in\B[\X\otimesw\Y,\Z]}$. So
$\phi=C_\theta(\Phi)={\J(\Phi)\sse\R(\J)}$. Thus, as this holds
for every ${\phi\in\bb[\X\x\X,\Z]}$, we get ${\bb[\X\x\X,\Z]\sse\R(\J)}$.

\vskip6pt\noi
(a)
By (a$_1$) and (a$_2$) we get ${\R(\J)=\bb[\X\x\Y,\Z]}$.

\vskip6pt\noi
(b)
Take an arbitrary ${\Phi\in\B[\X\otimesw\Y,\Z]}$. As we saw in $(*)$,
$$
\|\J(\Phi)(x,y)\|_{_\Z}\le\|\Phi\|_{_{_\ssB}}\|x\|_{_\X}\|y\|_{_\Y}
$$
for every ${(x,y)\in\X\x\Y}$. Then the bilinear map $\J(\Phi)$ is bounded
with norm
$$
\|\J(\Phi)\|_{_\ssb}\le \|\Phi\|_{_{_\ssB}}.
$$
Conversely, the unique ${\phi=C_\theta(\Phi)=\J(\Phi)}$ is such that
${\phi\in\bb[\X\x\Y,\Z]}$ because ${\R(\J)\sse\bb[\X\x\Y,\Z]}$ by (a$_1$)
and so, as we saw in $(**)$,
$$
\|\Phi(\digamma)\|_{_\Z}\le\|\phi\|_{_\ssb}\|\digamma\|_{_\we}
$$
for every ${\digamma\in\X\otimesw\Y}$. Therefore,
$$
\|\Phi\|_{_{_\ssB}}\le\|\phi\|_{_\ssb}=\|\J(\Phi)\|_{_\ssb}.
$$
Thus ${\|\J(\Phi)\|_{_\ssb}=\|\Phi\|_{_{_\ssB}}}$ for every
${\Phi\in\B[\X\otimesw\Y,\Z]}$. That is, $\J$ is an isometry.

\vskip6pt\noi
(c)
Then the linear transformation ${\J\!:\B[\X\otimesw\Y]\to\bb[\X\x\Y,\Z]}$
is a surjective isometry by (a) and (b), which means $\J$ is an isometric
isomorphism. Hence ${\B[\X\otimesw\Y,\Z]}$ and ${\bb[\X\x\Y,\Z]}$ are
isometrically isomorphic Banach spaces,
$$
\B[\X\otimesw\Y,\Z]\cong\bb[\X\x\Y,\Z].
$$
Consider the completion ${\X\hotimesw\Y}$ of ${\X\otimesw\Y}$. Since $\Z$ is
complete, it follows that
$$
\B[\X\hotimesw\Y,\Z]\cong\B[\X\otimesw\Y,\Z].
$$
So the stated result follows by transitivity. \qed

\section{Nuclear Operators}

\vskip6pt
Let ${\X,\Y,\V}$ be Banach spaces$.$ As a starting point, consider the
following expression involving the projective tensor product.
$$
\X^*\hotimes_\we\V^*\sse(\X\!\hotimesw\V)^*\!\cong\bb[\X\x\V,\FF]\cong
\B[\X,\V^*].
$$
The above inclusion comes from the definition of a reasonable crossnorm, the
first isometric isomorphism is a particular case of the {\it universal mapping
principle}\/ for the projective norm as in Theorem 3.2, and the second one is
the classical identification of bounded bilinear forms with bounded linear
transformations (see, e.g., \cite[Section 1.4, p.6]{DF})$.$ For the particular
case where $\Y$ is isometrically isomorphic to the dual of some Banach space
$\V$ (for instance, if is $\Y$ is reflexive), then we get
${\X^*\hotimes_\we\Y\sse\B[\X,\Y]}$, where the inclusion means algebraic
embedding$.$ On the other hand, this not only holds in general (no restriction
to Banach spaces being the dual of some Banach space) but is strengthened to
an isometric embedding for the injective norm instead$:$ {\it the injective
tensor product\/ ${\X^*\hotimesv\Y}$ is isometrically embedded in the Banach
space\/ $\B[\X,\Y]$, and so it is viewed as a subspace of}\/ $\B[\X,\Y]$:
$$
\X^*\hotimesv\Y\sse\B[\X,\Y]
$$
(see, e.g., \cite[Proposition 1.1.5]{DFS})$.$ This fails in general for the
projective norm, and the theorem below shows how far one can get along this
line in the general case.

\theorem{4.1}
{\it There is a natural transformation ${K\!:\!\X^*\hotimesw\Y\to\B[\X,\Y]}$
such that}\/
\begin{description}
\item{$\kern-12pt$(a)$\:$}
$K$ {\it is a linear contraction with}\/ ${\|K\|=1}$,
\vskip4pt
\item{$\kern-12pt$(b)$\:$}
{\it the range $\R(K)$ of $K$ is characterized as follows$\,:\!$
${T\kern-1pt\in\B[\X,\Y]}$ lies in\/ $\R(K)$ if and only if there are\/
$\X^*\!$-valued and\/ $\Y$-valued sequences\/ $\{f_k\}$ and\/ $\{y_k\}$ for
which the real sequence\/ $\{\|f_k\|\kern1pt\|y_k\|\}$ is summable and}
\vskip5pt\noi
$$
Tx={\sum}_k f_k(x)\kern1pt y_k
\quad\hbox{\it for every}\quad
x\in\X,
$$
{\it and}
\vskip4pt
\item{$\kern-12pt$(c)$\;$}
{\it for each ${T\!\in\kern-1pt\R(K)}$,
\vskip5pt\noi
$$
\|T\|_N
=\inf_{\{f_k\},\{y_k\},\,T=\Sigma_k f_k(\cdot)\,y_k}
{\sum}_k\|f_k\|\kern1pt\|y_k\|
$$
defines a norm on the linear space $\R(K)$ for which}\/
$\,{\|T\|\le\|T\|_N}$.
\end{description}

\proof
Let $\X^*$ be the dual of $\X.$ Take an arbitrary
$\digamma\!={\sum_kf_k\otimes y_k\in\!\X^*\hotimesw\Y}$ so that
$\|\digamma\|_\we=\inf{\sum_k\|f_k\|\kern1pt\|y_k\|}$ by Theorem 3.1$.$
Associated with $\digamma$ consider the natural transformation
${\Psi_\digamma\!:\X\to\Y}$ given by
\vskip6pt\noi
$$
\Psi_{\digamma}x={\sum}_kf_k(x)\,y_k
\quad\hbox{for every}\quad
x\in\X,
$$
which is linear and does not depend on the representation
${\sum_kf_k\otimes y_k}$ of $\digamma.$ Also, $\Psi_\digamma$ is bounded$.$ In
fact $\|\Psi_{\digamma}x\|\le{{\sum}_k\|f_k\|\kern1pt\|y_k\|\kern1pt\|x\|}.$
So ${\|\Psi_{\digamma}x\|\le \|\digamma\|_\we\|x\|}$, for every
${x\in\X}.$ Thus $\Psi_{\digamma}$ lies in ${\B[\X,\Y]}$ for each
${\digamma\in\X^*\hotimesw\Y}.$ This defines a transformation
$$
K\!:\X^*\hotimesw\Y\to\B[\X,\Y]
\quad\hbox{such that}\quad
K(\digamma)=\Psi_\digamma
\quad\hbox{for every}\quad
\digamma\in\X^*\hotimesw\Y.
$$

\vskip0pt\noi
(a)
$K$ is clearly linear$.$ It is contraction as well$.$ In fact, as
${\|\Psi_{\digamma}x\|\le \|\digamma\|_\we\|x\|}$ for every
${x\kern-1pt\in\kern-1pt\X}$,
$\,\|K(\digamma)\|={\|\Psi_{\digamma}\|=\sup_{\|x\|=1}\|\Psi_{\digamma}x\|
\le\|\digamma\|_\we}$ for every ${\digamma\!\in\kern-1pt\X^*\hotimesw\Y}.$
So ${K\in\B\big[\X^*\hotimesw\Y,\B[\X,\Y]\big]}$ with
${\|K\|=\sup_{\|\digamma\|_\we=1}\|K(\digamma)\|\le\kern-1pt1}.$ Reversely,
for ${f\otimes y}$ in ${\X^*\hotimesw\Y}$ with $\|f\|={\|y\|=1}$ we get
${\|\Psi_{f\otimes\,y}x\|}={|f(x)|\kern1pt\|y\|}={|f(x)|}$, and hence
${\|K(f\otimes y)\|}={\|\Psi_{f\otimes\,y}\|}={\sup_{\|x\|=1}|f(x)|}=\|f\|=1
={\|f\|\kern1pt\|y\|}={\|f\otimes y\|_\we}.$ Then with
${\digamma\!=f\otimes y}$ in ${\X^*\otimesw\Y}$ we get ${\|\digamma\|_\we=1}$
and $\|K(\digamma)\|=\|\digamma\|_\we.$ Thus ${\|K\|\ge1}$.

\vskip6pt\noi
(b)
By definition, ${\K(f\otimes y)=T\kern-1pt\in\B[\X,\Y]}$ if an only if
${Tx=\Psi_{f\otimes y}(x)=f(x)y}$ for every ${x\in\X}.$ Then for an
arbitrary ${\digamma\!=\sum_kf_k\otimes y_k\in\X^*\hotimesw\Y}$,
$$
\K(\digamma)=T\kern-1pt\in\B[\X,\Y]
\quad\hbox{if an only if}\quad
Tx=\Psi_{\digamma}(x)={\sum}_kf_k(x)\kern1pt y_k
$$
for every ${x\kern-1pt\in\kern-1pt\X}$, since
${K(\digamma)\kern-1pt=\!\sum_k\kern-1ptK(f_k\otimes y_k)}$ because $K$ is
linear and bounded$.$ This characterizes the range $\R(K)$ of $K$.

\vskip4pt\noi
(c)
As for ${\|\cdot\|_N}$ being a norm on $\R(K)$, we verify the triangle
inequality only$.$ Take ${T,T'\in\R(K)}$ so that
${T+T'}
={\sum_kf_k(\cdot)y_k}+{\sum_kf'_k(\cdot)y'_k}
={\sum_kf''_k(\cdot)y''_k}.$
Thus
${\|T+T'\|_N}\!
={\inf_{\,f_k,y_k,f'_k,y'_k}
(\sum_k\|f_k\|\kern1pt\|y_k\|+\|f'_k\|\kern1pt\|y'_k\|)}
={\inf_{\,f''_k,y''_k}\sum_k\|f''_k\|\kern1pt\|y''_k\|}
\le{\inf_{\,f_k,y_k}\sum_k\|f_k\|\kern1pt\|y_k\|}
\kern1pt+\kern1pt{\inf_{\,f'_k,y'_k}\sum_k\|f'_k\|\kern1pt\|y'_k\|}
={\|T\|_N+\|T'\|_N}.$
Now for the norm inequality$:$ if ${T\!\in\kern-1pt\R(K)}$, then there
are $\{f_k\}$ and $\{y_k\}$ with $\{\|f_k\|\kern1pt\|y_k\|\}$ summable such
that $T(x)={\sum_kf_k(x)y_k}$, and so $\|Tx\|={\|{\sum}_kf_k(x)y_k\|}$,
for every ${x\in\X}.$ Thus
$$
\|T\|
=\sup_{\|x\|=1}\|Tx\|
\le\sup_{\|x\|=1}{\sum}_k\|f_k\|\kern1pt\|y_k\|\kern1pt\|x\|
={\sum}_k\|f_k\|\kern1pt\|g_k\|
$$
so that ${\|T\|\le\inf\sum_k\|f_k\|\kern1pt\|g_k\|=\|T\|_N}$. \qed

\vskip9pt
A transformation ${T\kern-1pt\in\B[\X,\Y]}$ is {\it nuclear}\/ if it lies in
the range $\R(K)$ of such a natural contraction with unit norm
${K\!\in\B[\X^*\hotimesw\Y,\,\B[\X,\Y]\kern1pt]}$ defined in Theorem 4.1, and
the range $\R(K)$ is called {\it the linear space of nuclear
transformations}\/$.$
\hbox{Thus set}
$$
\B_N[\X,\Y]=\R(K)\sse\B[\X,\Y].
\quad\hbox{In particular,}\quad
\B_N[\X]=\B_N[\X,\X]\sse\B[\X].
$$
Theorem 4.1 prompts the following redefinition:
\vskip6pt\noi
{\narrower
A linear transformation ${T\kern-1pt\in\B[\X,\Y]}$ between Banach spaces
$\X$\and $\Y$ is {\it nuclear}\/ (i.e.,
${T\kern-1pt\in\kern-1pt\B_N[\X,\Y]}\kern.5pt)$ if there are $\X^*\!$-valued
and $\Y$-valued sequences $\{f_k\}$ and $\{y_k\}$ such that
${\sum_k\|f_k\|\kern1pt\|y_k\|<\infty}$ and $Tx={\sum_k f_k(x)\kern1pt y_k}$
for every ${x\in\X}$.
\vskip6pt}
\vskip0pt\noi
The expression $Tx={\sum_k f_k(x)\kern1pt y_k}$ is a {\it nuclear
representation}\/ of $T$ and
$$
\|T\|_N=\inf{\sum}_k\|f_k\|\kern1pt\|y_k\|
\quad\hbox{for every}\quad
T\kern-1pt\in\B_N[\X,\Y]
$$
(where the infimum is taken over all nuclear representations of $T$) defines
a norm on the linear space ${\B_N[\X,\Y]}$, the {\it nuclear norm}\/, such
that ${\|T\|\le\|T\|_N}.$

\vskip6pt
The linear contraction $K$ of norm one in Theorem 4.1 is not necessarily
injective (thus it may not be an isometry) and $\R(K)$ is identified
with the quotient space of ${\X^*\hotimesw\Y}$ modulo $\N(K)$, that is,
${\B_N[\X,\Y]}\cong{(\X^*\hotimesw\Y)/\N(K)}$ (see, e.g., \cite[p.41]{Rya}),
and therefore ${\B_N[\X,\Y]}$ is a Banach space$.$ However, if one of $\X^*\!$
or $\Y$ is a Banach space with the approximation property, then ${\N(K)=\0}$
(see, e.g.,\cite[Proposition 4.6]{Rya}), and so in this case
${\B_N[\X,\Y]}\cong{\X^*\hotimesw\Y}.$ In particular,
$$
\X\;\;\hbox{is a Hilbert space}
\;\limply\,
\B_N[\X]\cong\X^*\hotimesw\X\cong\X\hotimesw\X^*.
$$
Moreover, $\B_N[\X]$ is a two-sided ideal of the Banach algebra $\B[\X]$ by
the next result.

\corollary{4.2}
{\it Let\/ ${\X,\Y,\V,\W}$ be Banach spaces$.$ If\/
${T\kern-1pt\in\B_N[\X,\Y]}$,\/ ${R\in\B[\V,\X]}$ and\/ ${L\in\B[\Y,\W]}$,
then ${L\kern1ptTR\in\B_N[\V,\W]}$ and}\/
${\|L\kern1ptTR\|_N\le\|L\|\kern1pt\|T\|_N\|R\|}$.

\proof
This is a straightforward consequence of the definition of nuclear
transformation: as ${\sum_k\|f_k\|\kern1pt\|y_k\|<\infty}$ and
$Tx={\sum_kf_k(x)\kern1pt y_k}$, set $g_k={f_kR\in\V^*\kern-1pt}$ and
$w_k={Ly_k\in\W}$. So
${\sum_k\|g_k\|\kern1pt\|y_k\|
\le\sum_k\|f_k\|\kern1pt\|y_k\|\kern1pt\|R\|<\infty}$
and ${TRx=\sum_kg_k(x)\kern1pt y_k}.$ Also
${\sum_k\|f_k\|\kern1pt\|w_k\|}
\le{\sum_k\|f_k\|\kern1pt\|w_k\|\kern1pt\|L\|<\infty}$
and $LTx={\sum_kf_k(x)\kern1pt w_k}$. \qed

\vskip9pt
{From} now on let $\X$ be a Hilbert space$.$ Denote the inner product in $\X$
by ${\<\cdot\,;\cdot\>}$, and take the Banach algebra $\B[\X]$ of all
operators on $\X.$ Consider the Fourier Series and the Riesz Representation
Theorems (see, e.g., \cite[Theorems 5.48 and 5.62]{EOT})$.$ A functional $f$
lies in $\X^*$ if and only if there is a unique $z$ in $\X,$ called the Riesz
representation of $f$, such that ${f(x)=\<x\,;z\>}$ for every ${x\in\X}$ and
$\|z\|=\|f\|.$ Thus on a Hilbert space the previous redefinition can be
rewritten as follows:
\vskip6pt\noi
{\narrower
An operator ${T\kern-1pt\in\B[\X]}$ on a Hilbert space $\X$ is {\it nuclear}\/
(i.e., ${T\kern-1pt\in\B_N[\X]}\kern.5pt)$ if there are $\X$-valued sequences
$\{z_k\}$ and $\{y_k\}$ such that ${\sum_k\|z_k\|\kern1pt\|y_k\|<\infty}$ and
$Tx={\sum_k\<x\,;z_k\>\kern1pt y_k}$ for every ${x\in\X}$.
\vskip0pt}
\vskip6pt\noi
For each ${T\kern-1pt\in\B[\X]}$ let ${T^*\!\in\B[\X]}$ stand for its
Hilbert-space adjoint.

\corollary{4.3}
{\it Let\/ ${T\kern-1pt\in\B[\X]}$ be an operator on a Hilbert space}\/ $\X$.
\vskip6pt\noi
{\rm(a)}
{\it If there exist\/ $\X$-valued sequences\/ $\{z_k\}$ and\/ $\{y_k\}$ such
that}
$$
Tx={\sum}_k\<x\,;z_k\>y_k
\quad\hbox{\it for every}\quad
x\in\X,
$$
\vskip-4pt\noi
{\it then}
$$
T^*x={\sum}_k\<x\,;y_k\>z_k
\quad\hbox{\it for every}\quad
x\in\X.
$$
\vskip4pt\noi
{\rm(b)}
{\it $T$ is nuclear if an only if its adjoint\/ $T^*\!$ is nuclear and}\/
${\|T^*\|_N=\|T\|_N}$.

\proof
Let ${T\kern-1pt\in\B[\X]}$ be an operator on a Hilbert space $\X$.
\vskip4pt\noi
(a)
Suppose there are $\X$-valued sequences $\{z_k\}$ and $\{y_k\}$ such that
\vskip6pt\noi
$$
Tx={\sum}_k\<x\,;z_k\>\kern1pt y_k
\quad\hbox{for every}\quad
x\in\X.
$$
Then
${\<Tx\,;y\>}
={\big\<\sum_k\<x\,;z_k\>\kern1pt y_k;y\big\>}
={\sum_k\<x\,;z_k\>\kern1pt\<y_k;y\>}
={\sum_k\big\<x\,;\<y\,;y_k\>\kern1pt z_k\big\>}
={\big\<x\,;\sum_k\<y\,;y_k\>\kern1pt z_k\big\>}$
for every ${x,y\in\X}$, and therefore (by uniqueness of the adjoint)
\vskip6pt\noi
$$
T^*y={\sum}_k\<y\,;y_k\>\kern1pt z_k
\quad\hbox{for every}\quad
y\in\X.
$$
(b)
Immediate from item (a) and the definition of nuclear operator on $\X$. \qed

\section{Trace-Class Operators}

Let $\X$ be a Hilbert space$.$ A summary of the elementary expressions
required in this section goes as follows$.$ For each
${T\kern-1pt\in\kern-1pt\B[\X]}$ set ${|T|=(T^*T)^\half\kern-1pt\in\B[\X]}.$
Then
$$
\|Tx\|^2\!
=\<Tx\,;Tx\>
=\<T^*Tx\,;x\>
=\big\<|T|^2x\,;x\big\>
=\big\<|T|x\,;|T|x\big\>
=\big\|\kern1pt|T|x\big\|^{_{\scriptstyle2}}
$$
for every ${x\in\X}.$ So ${\|T\|=\big\|\kern1pt|T|\kern1pt\big\|}.$ Since
${\|T\|^2=\|T^*T\|=\|TT^*\|=\|T^*\|^2}$, we get
$$
\big\|\kern1pt|T|\kern1pt\big\|^{_{\scriptstyle2}}
=\|T\|^2
=\big\|\kern1pt|T|^2\big\|
=\big\|\kern1pt|T^*|^2\big\|
=\|T^*\|^2
=\big\|\kern1pt|T^*|\kern1pt\big\|^{_{\scriptstyle2}}.
$$
Moreover, since ${\|Q^\half\|^2=\|Q\|=\|Q^2\|^\half}$ for every nonnegative
operator ${Q\in\B[\X]}$,
$$
\big\|\kern1pt|T|^\half\big\|^{_{\scriptstyle2}}
=\big\|\kern1pt|T|\kern1pt\big\|
=\big\|\kern1pt|T|^2\big\|^{_{\scriptstyle\half}}.
$$
\vskip-4pt\noi
Hence for each ${x\in\X}$,
$$
\big\|\kern1pt|T|x\big\|^{_{\scriptstyle2}}
=\big\|\kern1pt|T|^\half|T|^\half x\big\|^{_{\scriptstyle2}}
\le\big\|\kern1pt|T|^\half\big\|^{_{\scriptstyle2}}
\big\|\kern1pt|T|^\half x\big\|^{_{\scriptstyle2}}
=\big\|\kern1pt|T|\kern1pt\big\|
\,\big\|\kern1pt|T|^\half x\big\|^{_{\scriptstyle2}}
=\|T\|\,\big\<|T|x\,;x\big\>.
$$
Now let $\{e_\gamma\}_{\gamma\in\Gamma}$ and $\{f_\gamma\}_{\gamma\in\Gamma}$
be arbitrary orthonormal bases for a Hilbert space $\X$, indexed by an
arbitrary nonempty index set $\Gamma$ (alternate notation:
$\{e_\gamma\}_\gamma$ or $\{e_\gamma\}$)$.$ By the Parseval identity (viz.,
${\|x\|^2=\sum_\gamma|\<x\,;e_\gamma\>|^2}$ for every ${x\in\X}$) we get
$$
{\sum}_{\gamma\in\Gamma}\|Te_\gamma\|^2
\!=\!{\sum}_{(\gamma,\delta)\in\Gamma^2}|\<Te_\gamma\,;f_\delta\>|^2
\!=\!{\sum}_{(\gamma,\delta)\in\Gamma^2}|\<T^*f_\delta\,;e_\gamma\>|^2
\!=\!{\sum}_{\delta\in\Gamma}\|T^*f_\delta\|^2
$$
whenever any of the families $\{\|Te_\gamma\|\}_\gamma$ or
$\{\|T^*f_\gamma\|\}_\gamma$ is square summable (i.e., if
${\sum_\gamma\|Te_\gamma\|^2<\infty}$ or
${\sum_\gamma\|T^*f_\gamma\|^2<\infty}$) for some orthonormal bases for $\X.$
Applying the above displayed identity to
${|T|^{_{\scriptstyle\half}}\!\in\B[\X]}$ instead of ${T\kern-1pt\in\B[\X]}$
we get
$$
{\sum}_\gamma\big\<|T|e_\gamma\,;e_\gamma\big\>
={\sum}_
\gamma\big\<|T|f_\gamma\,;f_\gamma\big\>.
$$
Thus if the family of positive numbers
$\big\{\big\<|T|e_\gamma\,;e_\gamma\big\>\big\}_{^{\scriptstyle\gamma}}
=\big\{\big\|\kern1pt|T|^\half e_\gamma\big\|^{_{\scriptstyle2}}
\big\}_{^{\scriptstyle\gamma}}$
is summable, then its sum does not depend on the orthonormal basis
$\{e_\gamma\}_\gamma$ for $\X$.
\vskip6pt\noi
{\narrower
An operator ${T\kern-1pt\in\B[\X]}$ on a Hilbert space $\X$ is
$\kern.5pt${\it trace-class}$\kern1pt$ if
${\sum_\gamma\!\big\<|T|e_\gamma;e_\gamma\big\>\kern-1pt<\kern-1pt\infty}$
for an arbitrary orthonormal basis $\{e_\gamma\}$ for $\X$.
\vskip6pt}
\vskip0pt\noi
Let $\B_1[\X]$ denote the subset of $\B[\X]$ consisting of all trace-class
operators$.$ Since
${\big\<|T|x\,;x\big\>=\big\|\kern1pt|T|^\half x\big\|^{_{\scriptstyle2}}}$
for every ${x\in\X}$, an operator $T$ lies in $\B_1[\X]$ if and only if
${\sum_\gamma\big\|\kern1pt|T|^\half e_\gamma\big\|^{_{\scriptstyle2}}<\infty}$
for an arbitrary orthonormal basis $\{e_\gamma\}.$ So for
${T\kern-1pt\in\B_1[\X]}$ set
$$
\|T\|_1
={\sum}_\gamma\big\<|T|e_\gamma\,;e_\gamma\big\>
={\sum}_\gamma\big\|\kern1pt|T|^\half e_\gamma\big\|^{_{\scriptstyle2}}.
$$
Clearly, ${\big|\kern1pt|T|\kern1pt\big|=|T|}.$ Then
${|T|^2=T^*T\kern-1pt\in\B[\X]}$ is trace-class (i.e., $|T|^2\!$ lies in
$\B_1[\X]\kern.5pt$) if and only if
${\sum_\gamma\big\|\kern1pt|T|e_\gamma\big\|^{_{\scriptstyle2}}\!<\infty}$;
equivalently, if ${\sum_\gamma\|Te_\gamma\|^2\!<\kern-1pt\infty}$.
\vskip6pt\noi
{\narrower
An operator ${T\!\in\kern-1pt\B[\X]}$ on a Hilbert space $\X$ is
{\it Hilbert--Schmidt}\/ if ${|T|^2\!\in\kern-1pt\B_1[\X]}$; equivalently if
${\sum_\gamma\|Te_\gamma\|^2\!<\kern-1pt\infty}$ for an arbitrary orthonormal
basis $\{e_\gamma\}$ for $\X$.
\vskip6pt}
\vskip0pt\noi
Let $\B_2[\X]$ denote the subset of $\B[\X]$ consisting of all
Hilbert--Schmidt operators$.$ If $T$ lies in $\B_2[\X]$, then set
$\|T\|^2_2={\sum_\gamma\|Te_\gamma\|^2}$ so that from what we have seen above,
$$
\|T\|_2
=\Big({\sum}_\gamma\|Te_\gamma\|^2\Big)^{_{\scriptstyle\half}}
=\Big({\sum}_\gamma\big\|\kern1pt|T|e_\gamma\big\|^{_{\scriptstyle2}}
\Big)^{_{\scriptstyle\half}}
=\big\|\kern1pt|T|^2\big\|^{_{\scriptstyle\half}}_1
=\|T^*T\|^\half_1,
$$
and so (recall: ${\big|\kern1pt|T|\kern1pt\big|=|T|}\kern.5pt$) we may infer:
$$
T\kern-1pt\in\B_1[\X]
\!\iff\!
|T|\in\B_1[\X]
\!\iff\!
|T|^\half\kern-1pt\in\B_2[\X]
\quad\hbox{and}\quad
\|T\|_1
=\big\|\kern1pt|T|\kern1pt\big\|_1
=\big\|\kern1pt|T|^\half\big\|^{_{\scriptstyle2}}_2,
$$
$$
T\kern-1pt\in\B_2[\X]
\!\iff\!
|T|\in\B_2[\X]
\!\iff\!
|T|^2\kern-1pt\in\B_1[\X]
\quad\hbox{and}\quad
\|T\|^2_2
=\big\|\kern1pt|T|\kern1pt\big\|^2_{_{\scriptstyle2}}
=\big\|\kern1pt|T|^2\big\|_1.
$$

\lemma{5.1}
{\it Let\/ $\X$ be a Hilbert space$.$ The following assertions hold true}$.$
\begin{description}
\item{$\kern-12pt$\rm(a)}
{\it $\,$The set\/ $\B_2[\X]$ is a linear space and\/
${\|\cdot\|_2\!:\B_2[\X]\to\RR}$ is a norm on}\/ $\B_2[\X].$
\vskip4pt
\item{$\kern-12pt$\rm(b)}
{\it $\,$If\/ ${S\in\B[\X]}$ and\/ ${T\kern-1pt\in\B_2 [\X]}$, then}\/
${\max\{\|ST\|_2,\|TS\|_2\}\le\|S\|\,\|T\|_2}.$
\vskip4pt
\item{$\kern-12pt$\rm(c)}
{\it $\;\B_2[\X]$ is an ideal of\/ the algebra}\/ $\B[\X].$
\vskip4pt
\item{$\kern-12pt$\rm(d)}
{\it $\;\B_2[\X]\sse\B_\infty[\X]\;$ and\/ $\;{\|T\|\le\|T\|_2}$ for every}\/
${T\kern-1pt\in\B_2[\X]}.$
\vskip4pt
\item{$\kern-12pt$\rm(e)}
{\it $\;{T^*\!\in\B_2[\X]}$ if and only if\/ ${T\kern-1pt\in\B_2[\X]}\;$
and}\/ $\;{\|T\|_2=\|T^*\|_2}$.
\end{description}

\proof
Suppose $S$ and $T$ lie in $\B[\X]$ and consider the following assertion.

\claim{1}
If ${S,T\kern-1pt\in\B_2[\X]}$, then ${S+T\kern-1pt\in\B_2[\X]}$ and
${\|S+T\|_2\le\|S\|_2+\|T\|_2}$.

\vskip6pt\noi
{\it Proof of Claim 1}$.$
Let the index $\gamma$ run over an arbitrary index set ${\Gamma\ne\void}$,
consider the Schwarz inequality on the Hilbert space $\ell^2$ over $\Gamma$,
and take ${S,T\kern-1pt\in\B_2[\X]}$ \hbox{so that}
${\sum_\gamma\|Se_\gamma\|\kern1pt\|Te_\gamma\|
\kern-1pt\le\kern-1pt
\big(\sum_\gamma\|Se_\gamma\|^2\big)^{_{\scriptstyle\half}}
\big(\sum_\gamma\|Te_\gamma\|^2\big)^{_{\scriptstyle\half}}
\!=\kern-1pt\|S\|_2\|T\|_2}$ if ${S,T\!\in\kern-1pt\B_2[\X]}.$ Thus
\begin{eqnarray*}
\|T+S\|^2_2
&\kern-6pt=\kern-6pt&
{\sum}_\gamma\|Se_\gamma+Te_\gamma\|^2
\le{\sum}_\gamma\big(\|Se_\gamma\|+\|Te_\gamma\|\big)^2                  \\
&\kern-6pt=\kern-6pt&
{\sum}_\gamma\|Se_\gamma\|^2+{\sum}_\gamma\|Te_\gamma\|^2
+2{\sum}_\gamma\|Se_\gamma\|\,\|Te_\gamma\|                              \\
&\kern-6pt\le\kern-6pt&
\|S\|^2_2+\|T\|^2_2+2\|S\|_2\|T\|_2=(\|S\|_2+\|T\|_2)^2. \qed
\end{eqnarray*}
\vskip-2pt

\vskip9pt\noi
Since homogeneity, nonnegativity and positivity for
${\|\cdot\|_2\!:\B_2[\X]\to\RR}$ is readily verified, Claim 1 is enough to
ensure that
$$
\hbox{$\B_2[\X]$ is a linear space and ${\|\cdot\|_2}$ is a norm on it.}
$$
Also, if ${S\in\B[\X]}$ and ${T\kern-1pt\in\B_2 [\X]}$, then
${\|S\kern1pt T\|^2_2}\kern-1pt
={\sum_\gamma\|S\kern1pt Te_\gamma\|^2}\kern-1pt
\le{\|S\|^2\sum_\gamma\|Te_\gamma\|^2}\kern-1pt
={\|S\|^2\,\|T\|^2_2}.$
Similarly,
${\|TS\|^2_2}
={\sum_\gamma\|TS e_\gamma\|^2}
={\sum_\gamma\|(TS)^*e_\gamma\|^2}
={\sum_\gamma\|S^*T^*e_\gamma\|^2}
\le{\|S^*\|^2\sum_\gamma\|T^*e_\gamma\|^2}
={\|S\|^2\sum_\gamma\|Te_\gamma\|^2}={\|S\|\,\|T\|^2_2}.$
Therefore,
$$
\max\{\|ST\|_2,\|TS\|_2\}\le\|S\|\,\|T\|_2
$$
\vskip-6pt\noi
and
$$
\hbox{$\B_2[\X]$ is an ideal (i.e., a two-sided ideal) of $\B[\X]$.}
$$
\vskip4pt\noi
Now take an arbitrary ${T\!\in\kern-1pt\B_2[\X]}$ so that
${\sum_{\gamma\in\Gamma}\|T^*e_\gamma\|^2}\!
={\sum_{\gamma\in\Gamma}\|Te_\gamma\|^2}\!<\kern-1pt\infty$,
and take any integer
${n\kern-1pt\ge\kern-1pt1}.$ Thus there is a finite set
${N_n\kern-1pt\sse\kern-1pt\Gamma}$ such that
${\sum_{k\in N}\!\|T^*e_k\|^2\!<\kern-1pt\frac{1}{n}}$ for all finite sets
${N\kern-1pt\sse\kern-1pt\Gamma\\N_n}$ (Cauchy criterion for summable
families --- see, e.g.,\cite[Theorem 5.27]{EOT})$.$ So
${\sum_{\gamma\in\Gamma\\N_n}\!\|T^*e_\gamma\|^2<\frac{1}{n}}.$ Recall that
$Tx={\sum_{\gamma\in\Gamma}\<Tx\,;e_\gamma\>e_\gamma}$ for every
${x\kern-1pt\in\kern-1pt\X}$ (Fourier series expansion)$.$ Set
${T_nx=\kern-1pt\sum_{k\in N_n}\kern-1pt\<Tx\,;e_k\>e_k}$ for each
${x\kern-1pt\in\kern-1pt\X}$, which defines an operator $T_n$ in $\B_0[\X]$
because $N_n$ is finite$,$ Hence
${\|(T-T_n)x\|^2}={\sum_{\gamma\in\Gamma\\N_n}|\<Tx\,;e_\gamma\>|^2}
\le\big({\sum_{\gamma\in\Gamma\\N_n}\|T^*e_\gamma\|^2\big)\|x\|^2}$
for every ${x\kern-1pt\in\kern-1pt\X}.$ This implies that ${\|T_n-T\|\to0}$,
and so $T$ is the uniform limit of a sequence of finite-rank operators on a
Banach space, and therefore $T$ is compact (see, e.g.,
\cite[Corollary 4.55]{EOT}).
$$
\hbox{Every Hilbert--Schmidt operator is compact.}
$$
As we saw above,
${\sum_\gamma\|Te_\gamma\|^2=\sum_\gamma\|T^*e_\gamma\|^2}$ and
${\|Te\|\le\|T\|_2}$ if ${\|e\|=1}$, and so
$$
\hbox{ ${T^*\!\in\B_2[\X]}$ if and only if ${T\kern-1pt\in\B_2[\X]}\;$ and
$\;{\|T\|\le\|T\|_2=\|T^*\|_2}$}.                                 \eqno{\qed}
$$
\vskip-2pt

\vskip9pt
The norm ${\|\cdot\|_2}$ on the linear space $\B_2[\X]$ of all Hilbert--Schmidt
operators is referred to as the {\it Hilbert--Schmidt norm}\/.

\theorem{5.2}
{\it Let\/ $\X$ be a Hilbert space$.$ The following assertions hold true}$.$
\begin{description}
\item{$\kern-12pt$\rm(a)}
{\it $\,$The set\/ $\B_1[\X]$ is a linear space and\/
${\|\cdot\|_1\!:\B_1[\X]\to\RR}$ is a norm on}\/ $\B_1[\X].$
\vskip4pt
\item{$\kern-12pt$\rm(b)}
{\it $\;{T\kern-1pt\in\B_1[\X]}$ if and only if\/ ${T=A\kern1pt B}$ for
some}\/ ${A,B\in\B_2[\X]}.$
\vskip4pt
\item{$\kern-12pt$\rm(c)}
{\it $\;{\B_1[\X]\sse\B_2[\X]}\,$ and\/ $\,{\|T\|^2_2\le\|T\|\,\|T\|_1}$ so
that\/ ${\|T\|_2\le\|T\|_1}$ for}\/ ${T\kern-1pt\in\B_1[\X]}.$
\vskip4pt
\item{$\kern-12pt$\rm(d)}
{\it $\;\B_1[\X]$ is an ideal of the algebra}\/ $\B[\X].$
\vskip4pt
\item{$\kern-12pt$\rm(e)}
{\it $\,$If\/ ${S\in\B[\X]}$ and\/ ${T\!\in\B_1[\X]}$, then}\/
${\max\{\|ST\|_1,\|TS\|_1\}\le\|S\|\,\|T\|_1}.$
\vskip4pt
\item{$\kern-12pt$\rm(f)}
{\it $\;{T^*\!\in\B_1[\X]}$ if and only if\/ ${T\kern-1pt\in\B_1[\X]}\;$
and}\/ $\;{\|T\|_1=\|T^*\|_1}.$\vskip6pt\noi
\item{$\kern-12pt$\rm(g)}
$\;{\B_0[\X]\sse\B_1[\X]}$.
\end{description}

\proof
Suppose $S$ and $T$ lie in $\B[\X]$ and consider the following assertion.

\claim{2}
If ${S,T\kern-1pt\in\B_1[\X]}$, then ${S+T\kern-1pt\in\B_1[\X]}$ and
${\|S+T\|_1\le\|S\|_1+\|T\|_1}$.

\vskip6pt\noi
{\it Proof of Claim 2}$.$
Consider the polar decompositions ${T\!+\kern-1ptS=W|T\kern-1pt+S|}$,
$\,{T\kern-1pt=W_1|T|}$ and $S={W_2|S|}$, where
${W,\kern-1ptW_1,\kern-1ptW_2}$ are partial isometries in ${\B[\X]}$, so that
${|T\!+\kern-1ptS|}={W^*(T\!+\kern-1ptS)}$, $\,|T|={W_1^*T}$, and
$|S|={W_2^*T}.$ Hence (Schwartz inequality on $\ell^2$ over $\Gamma$)
\begin{eqnarray*}
\|T+S\|_1
&\kern-6pt=\kern-6pt&
{\sum}_\gamma\big\<|T+S|e_\gamma;e_\gamma\big\>
\le{\sum}_\gamma|\<Te_\gamma;W^*e_\gamma\>|
\,+{\sum}_\gamma|\<Se_\gamma;W^*e_\gamma\>|                               \\
&\kern-6pt=\kern-6pt&
{\sum}_\gamma
\big|\big\<|T|^\half e_\gamma\,;|T|^\half W_1^*We_\gamma\big\>\big|
\,+{\sum}_\gamma
\big|\big\<|S|^\half e_\gamma\,;|S|^\half W_2^*We_\gamma\big\>\big|       \\
&\kern-6pt\le\kern-6pt&
{\sum}_\gamma
\big\|\kern1pt|T|^\half e_\gamma\big\|
\,\big\|\kern1pt|T|^\half W_1^*We_\gamma\big\|
\,+{\sum}_\gamma
\big\|\kern1pt|S|^\half e_\gamma\big\|
\,\big\|\kern1pt|S|^\half W_2^*We_\gamma\big\|                            \\
&\kern-6pt\le\kern-6pt&
\Big({\sum}_\gamma
\big\|\kern1pt|T|^\half e_\gamma\big\|^{_{\scriptstyle2}}
\Big)^{_{\scriptstyle\half}}
\Big({\sum}_\gamma
\big\|\kern1pt|T|^\half W_1^*We_\gamma\big\|^{_{\scriptstyle2}}
\Big)^{_{\scriptstyle\half}}                                              \\
&\kern-6pt+\kern-6pt&
\Big({\sum}_\gamma
\big\|\kern1pt|S|^\half e_\gamma\big\|^{_{\scriptstyle2}}
\Big)^{_{\scriptstyle\half}}
\Big({\sum}_\gamma
\big\|\kern1pt|S|^\half W_2^*We_\gamma\big\|^{_{\scriptstyle2}}
\big)^{_{\scriptstyle\half}}                                              \\
&\kern-6pt=\kern-6pt&
\big\|\kern1pt|T|^\half\big\|_2\;\big\|\kern1pt|T|^\half W^*_1W\big\|_2
+\big\|\kern1pt|S|^\half\big\|_2\;\big\|\kern1pt|S|^\half W^*_1W\big\|_2  \\
&\kern-6pt\le\kern-6pt&
\big\|\kern1pt|T|^\half\big\|_2^{_{\scriptstyle2}}\;\|W_1^*W\|
+\big\|\kern1pt|S|^\half\big\|_2^{_{\scriptstyle2}}\;\|W_2^*W\|
\;\le\;\|T\|_1+\|S\|_1.
\end{eqnarray*}
(Recall that
${T\kern-1pt\in\B_1[\X]}
\limply{|T|^\half\kern-1pt\in\B_2[\X]}
\limply{|T|^\half A\in\B_2[\X]}$,
$\,{\sum_\gamma\big\|\kern1pt|T|^\half e_\gamma\big\|^{_{\scriptstyle2}}}\!
=\big\|\kern1pt|T|^\half\big\|^{_{\scriptstyle2}}_2\!
=\|T\|_1$,
${\big(\sum_\gamma\!
\big\|\kern1pt|T|^\half Ae_\gamma\big\|^{_{\scriptstyle2}}
\big)^{_{\scriptstyle\half}}\!
=\big\|\kern1pt|T|^\half A\big\|_2\!
\le\big\|\kern1pt|T|^\half\big\|_2\|A\|}$
by Lemma 5.1(b), and $\|W\|=\|W_1\|=\|W_2\|=1$ since these are
partial isometries)$.\qed$

\vskip6pt\noi
Claim 2 is enough to ensure that
$$
\hbox{$\B_1[\X]$ is a linear space and ${\|\cdot\|_1}$ is a norm on it,}
$$
since ${\|\cdot\|_1}$ is trivially homogeneous (so every multiple of an
operator in $\B_1[\X]$ lies in $\B_1[\X]$), and nonnegativeness and
positiveness for ${\|\cdot\|_1}$ are readily verified$.$ Consider again the
polar decomposition ${T\kern-1pt=W_1|T|=W_1|T|^\half|T|^\half}.$ If
${T\kern-1pt\in\B_1[\X]}$, then $|T|^\half$ lies in $\B_2[\X]$ and so
${W_1|T|^\half}$ lies in $\B_2[\X]$ as well according to Lemma 5.1(c)$.$
Conversely, suppose ${T\kern-1pt=A\kern1ptB}$ with ${A,B\in\B_2[\X]}.$ Thus
$T$ lies in ${\B_2[\X]}$ by Lemma 5.1(c) again$.$ Since $|T|={W_1^*T}$ we get
$|T|={W^*\!A\kern1ptB}$ with ${A^*W\!\in\B_2[\X]}$ (according to
Lemma 5.1(c) once again)$.$ Therefore
$$
{\sum}_\gamma\big\<|T|e_\gamma\kern1pt;e_\gamma\big\>
\le{\sum}_\gamma\|Be_\gamma\|\|A^*We_\gamma\|
\le\Big({\sum}_\gamma\|Be_\gamma\|^2\Big)^{_{\scriptstyle\half}}
\Big({\sum}_\gamma\|A^*We_\gamma\|^2\Big)^{_{\scriptstyle\half}}
$$
(by the Schwarz inequality on both Hilbert spaces $\X$ and $\ell^2$ over
$\Gamma$, as we did before)$.$ Hence ${T\kern-1pt\in\B_1[\X]}$ with
${\|T\|_1\le\|B\|_2\|A^*W\|_2\le\|B\|_2\|A\|_2}.$ Summing up:
$$
\hbox{$T\kern-1pt\in\B_1[\X]\!\iff\!T\!=A\kern1pt B$ for some
$A,B\in\B_2[\X]\;$ and $\;\|T\|_1\le\|A\|_2\|B\|_2$}.
$$
Also, since the product ${A\kern1ptB}$ lies in ${\B_2[\X]}$ by Lemma 5.1(c),
$\,{\B_1[\X]\sse\B_2[\X]}$:
$$
\hbox{Every trace-class operator is Hilbert--Schmidt.}
$$
Moreover, since
${\big\|\kern1pt|T|x\big\|^{_{\scriptstyle2}}
\!\le\|T\|\,\big\|\kern1pt|T|^\half x\big\|^{_{\scriptstyle2}}}\kern-1pt$
for every ${x\in\kern-1pt\X}$ and ${\|T\|\le\|T\|_2}$
(by Lem\-ma 5.1(d)$\kern.5pt$), we get
${\|T\|^2_2
=\sum_\gamma\big\|\kern1pt|T|e_\gamma\big\|^{_{\scriptstyle2}}\kern-1pt
\le\|T\|\sum_\gamma\big\|\kern1pt|T|^\half e_\gamma\big\|^{_{\scriptstyle2}}
\kern-1pt=\|T\|\,\|T\|_1}.$ So
$$
\|T\|^2_2\le\|T\|\|T\|_1
\;\limply\;
\|T\|_2\le\|T\|_1
\quad\hbox{for every}\quad
T\kern-1pt\in\B_1[\X].
$$
Hence if $T$ lies in ${\B_1[\X]}$ and $S$ lies in ${\B[\X]}$, then
${S\kern1pt T=(SA)B}$ and ${TS=A(BS)}$ for some ${A,B\in\B_2[\X]}.$ Thus
${S\kern1ptT=CB}$ and ${TS=AD}$ with both ${C=SA}$ and ${D=BS}$ in
${\B_2[\X]}$ according to Proposition A1(c)$.$ So ${S\kern1ptT}$ and
${TS}$ lie in ${\B_1[\X]}.$ Therefore
$$
\hbox{$\B_1[\X]$ is an ideal (i.e., a two-sided ideal) of $\B[\X]$.}
$$
\vskip-4pt

\claim{3}
Let $\{e_\gamma\}$ be any orthonormal basis for $\X.$ If ${T\in\B_1[\X]}$ and
${S\in\B[\X]}$, then
\begin{description}
\item{$\kern-9pt$\rm(i)}
$\;{\sum_\gamma|\<Te_\gamma\,;e_\gamma\>|\le\|T\|_1}\;\;$ and
$\;\;{\sum_\gamma\<Te_\gamma\,;e_\gamma\>}$
does not depend on $\{e_\gamma\}$,
\vskip4pt
\item{$\kern-11pt$\rm(ii)}
$\;\sum_\gamma\<TSe_\gamma\,;e_\gamma\>
=\sum_\gamma\<S\kern1pt Te_\gamma\,;e_\gamma\>$,
\vskip4pt
\item{$\kern-12pt$\rm(iii)}
$\;\big|\sum_\gamma\<S|T|e_\gamma\,;e_\gamma\>\big|
=\big|\sum_\gamma\big\<|T|Se_\gamma\,;e_\gamma\big\>\big|
\le\|S\|\,\|T\|_1$.
\end{description}

\vskip6pt\noi
{\it Proof of Claim 3}$.$
Let ${T\kern-1pt=W_1|T|}$ be the polar decomposition of
${T\kern-1pt\in\B_1[\X]}.$ Since ${|T|^\half\in\B_2[\X]}$ with
${\big\|\kern1pt|T|^\half\big\|_2=\|T\|^\half_1}$ and ${\|W^*_1\|=1}$ we get
by Lemma 5.1(b),
\begin{eqnarray*}
{({\rm i}_1)\kern17pt}
{\sum}_\gamma|\<Te_\gamma\,;e_\gamma\>|
&\kern-6pt=\kern-6pt&
{\sum}_\gamma
\big|\big\<|T|^\half e_\gamma\,;|T|^\half W_1^*e_\gamma\big\>\big|
\le{\sum}_\gamma\big\|\kern1pt|T|^\half e_\gamma\big\|\,
\big\|\kern1pt|T|^\half W_1^*e_\gamma\big\|                              \\
&\kern-6pt\le\kern-6pt&
\Big({\sum}_\gamma
\big\|\kern1pt|T|^\half e_\gamma\big\|^{_{\scriptstyle2}}
\Big)^{_{\scriptstyle\half}}
\Big({\sum}_\gamma
\big\|\kern1pt|T|^\half W_1^*e_\gamma\big\|^{_{\scriptstyle2}}
\Big)^{_{\scriptstyle\half}}                                             \\
&\kern-6pt\le\kern-6pt&
\|T\|_1^\half\,\big\|\kern1pt|T|^\half W_1^*\big\|_2
\le\|T\|_1^\half\,\big\|\kern1pt|T|^\half\big\|_2\,\|W_1^*\|
\le\|T\|_1
\end{eqnarray*}
(by the Schwarz inequality on both Hilbert spaces $\X$ and $\ell^2$ over
$\Gamma$ again)$.$ Thus $\{\<Te_\gamma\,;e_\gamma\>\}$ is a summable family,
and by taking an arbitrary orthonormal basis $\{f_\gamma\}$ for the Hilbert
space $\X$, and applying the Fourier Series Theorem, we get
$$
{({\rm i}_2)\kern18pt}
{\sum}_\gamma\<Te_\gamma\,;e_\gamma\>
={\sum}_{\gamma,\delta}\<Te_\gamma\,;f_\delta\>\,\<f_\delta\,;e_\gamma\>,
={\sum}_\delta\<f_\delta\,;T^*f_\delta\>
={\sum}_\delta\<Tf_\delta\,;f_\delta\>.
$$
By (i$_1$) and (i$_2$) we get (i)$.$ As ${S\kern1ptT}$ and ${TS}$ lie in
${\B_1[\X]}$ by (d), it follows by (i) that the sums in (ii) exist and do not
depend on the orthonormal basis$.$ So let $\{f_\gamma\}$ be any orthonormal
basis for $\X$ and consider again the Fourier Series Theorem$.$ Thus
\begin{eqnarray*}
{\rm(ii)\kern55pt}
&&
{\kern-15pt\sum}_\gamma\<TSe_\gamma\,;e_\gamma\>
={\!\sum}_\gamma\<Se_\gamma\,;T^*e_\gamma\>
={\!\sum}_{\gamma,\delta}\<Se_\gamma\,;f_\delta\>\<Tf_\delta\,;e_\gamma\>  \\
&\kern-6pt=\kern-6pt&
\!{\sum}_{\gamma,\delta}\<Te_\gamma\,;f_\delta\>\<Sf_\delta\,;e_\gamma\>
={\!\sum}_\gamma\<Te_\gamma\,;S^*e_\gamma\>
={\!\sum}_\gamma\<S\kern1pt Te_\gamma\,;e_\gamma\>.
\end{eqnarray*}
Recall:
${T\kern-1pt\in\kern-1pt\B_1[\X]}
\!\!\iff\!\!{|T|\kern-1pt\in\kern-1pt\B_1[\X]}
\!\!\iff\!\!{|T|^\half\kern-1pt\in\kern-1pt\B_2[\X]}$
with
${\|T\|_1
\!=\big\|\kern1pt|T|\kern1pt\big\|_1
=\!\big\|\kern1pt|T|^\half\big\|^{_{\scriptstyle2}}_2}$
(so that ${S|T|}$ and ${|T|S}$ lie in ${\B_1[\X]}$ by (d)$\kern.5pt$), and if
${A\in\B_2[\X]}$ and ${B\in\B[\X]}$, then ${\|AB\|_2}\le{\|A\|_2\|B\|}$ by
Lemma 5.1(b)$.$ Therefore, before applying (ii), we get
\begin{eqnarray*}
{\rm(iii)\kern52pt}
\Big|{\sum}_\gamma\<S|T|e_\gamma\,;e_\gamma\>\Big|
&\kern-6pt=\kern-6pt&
\Big|{\sum}_\gamma\big\<|T|^\half e_\gamma\,;|T|^\half S^*e_\gamma\big\>
\Big|                                                                      \\
&\kern-6pt\le\kern-6pt&
\Big({\sum}_\gamma\big\|\kern1pt|T|^\half e_\gamma\big\|^{_{\scriptstyle2}}
\Big)^{_{\scriptstyle\half}}
\Big({\sum}_\gamma\big\|\kern1pt|T|^\half S^*e_\gamma\big\|^{_{\scriptstyle2}}
\Big)^{_{\scriptstyle\half}}                                               \\
&\kern-6pt=\kern-6pt&
\big\|\kern1pt|T|^\half\big\|_2\;\big\|\kern1pt|T|^\half S^*\big\|_2
\le\big\|\kern1pt|T|^\half\big\|_2\;\big\|\kern1pt|T|^\half
\big\|_2\;\|S^*\|                                                          \\
&\kern-6pt=\kern-6pt&
\|T\|_1\;\|S\|. \qed
\end{eqnarray*}
\vskip-2pt

\vskip6pt\noi
Now apply Claim 3(iii) to support the following argument$.$ Consider again the
polar decompositions ${T=W_1|T|}$, $\,{S\kern1ptT}={W_L|S\kern1ptT|}$, and
${TS=W_R|TS|}$, where $W_1$, $W_L,$ $W_R$ are partial isometries in $\B[\X]$
(with norm one as well as their adjoint), and so
${|S\kern1ptT|}={W_L^*S\kern1ptT}={W_L^*SW_1|T|}$ and
${|TS|}={W_R^*TS}={W_R^*W_1|T|S}.$ Since $T$ lies in $\B_1[\X]$, Claim 3(iii)
ensures
${\|S\kern1ptT\|_1}
={\sum_\gamma\big\<|S\kern1ptT|e_\gamma\,;e_\gamma\big\>}
={\sum_\gamma\<W_L^*SW_1|T|e_\gamma\,;e_\gamma\>}
\le{\|W_L^*SW_1\|\,\|T\|_1}
\le
{\|W_L^*\|\,\|S\|\,\|W_1\|\,\|T\|_1}
={\|S\|\,\|T\|_1}.$
Since ${W_R^*W_1|T|\in\B_1[\X]}$, Claim 3(iii) also ensures
${\|TS\|_1}
={\sum_\gamma\big\<|TS|e_\gamma\,;e_\gamma\big\>}
={\sum_\gamma\<W_R^*W_1|T|Se_\gamma\,;e_\gamma\>}
={\sum_\gamma\<SW_R^*W_1|T|e_\gamma\,;e_\gamma\>}\kern-1pt
\le{\|SW_R^*W_1\|\,\|T\|_1}\kern-1pt
\le{\|S\|\,\|W_R^*\|\,\|W_1\,\|T\|_1}\kern-1pt
={\|S\|\,\|T\|_1}.\kern-1pt$ Thus
$$
\max\{\|ST\|_1,\|TS\|_1\}\le\|S\|\,\|T\|_1.
$$
If ${T\kern-1pt\in\B_1[\X]}$, then ${T\kern-1pt=\kern-1ptA\kern1ptB}$ with
${A,B\in\B_2[\X]}$ by (b) (and so ${A^*\!,B^*\kern-1pt\in\B_2[\X]}$ by
Lemma 5.1(e)$\kern.5pt).$ Then, using (b) again, ${T^*\!=B^*A^*\in\B_1[\X]}.$
Dually, if ${T^*\in\B_1[\X]}$, then ${T=T^{**}\in\B_1]\X]}.$ Now by taking the
polar decompositions ${T\kern-1pt=W_1|T|}$ and ${T^*\!=W'_1|T^*|}$ we get
$|T^*|={W^{\prime*}_1T^*\kern-1pt}={W^{\prime*}_1|T|W^*_1}.$ Therefore
$\|T^*\|_1
=\big\|\kern1pt|T^*|\kern1pt\big\|_1
={\|W^{\prime*}_1|T|W^*_1\|_1}
\le{\|W^{\prime*}_1\|\,\big\|\kern1pt|T|\kern1pt\big\|_1\,\|W^*_1\|_1}
=\big\|\kern1pt|T|\kern1pt\big\|_1
=\|T\|_1$
by (e) (proved above)$.$ Dually,
${\|T\|_1=\|T^{**}\|_1\kern-1pt\le\|T^*\|_1}.$ Thus
$$
T\kern-1pt\in\B_1[\X]\!\iff\!T^*\!\in\B_1[\X]
\quad\hbox{and}\quad
\|T^*\|_1\kern-1pt=\|T\|_1.
$$
Finally, recall that ${\R(T)^-\!=\N(T^*)^\perp}.$ If $\dim\R(T)$ is finite,
then so is $\dim\N(T^*)^\perp.$ Let $\{e_\delta\}$ be an orthonormal basis for
$\N(T^*)$ and let $\{e_k\}$ be a finite orthonormal basis for
$\N(T^*)^\perp\kern-1pt.$ As ${\X=\N(T^*)+\N(T^*)^\perp}$, then
${\{e_\gamma\}=\{e_\delta\}\cup\{e_k\}}$ is an orthonor\-mal basis for $\X.$
Now either ${T^*e_\gamma=0}$ or ${T^*e_\gamma=T^*e_k}.$ Therefore
${{\sum}_\gamma\big\<|T^*|e_\gamma;e_\gamma\big\>}
={{\sum}_k\big\<|T^*|e_k;e_k\big\>}
<\infty.$
Thus ${T^*\!\in\B_1[\X]}.$ So ${T\kern-1pt\in\B_1[\X]}$ by (f)$.$ Hence
${\B_0[\X]\sse\B_1[\X]}$:
$$
\hbox{Every finite-rank operator is trace-class}. \eqno{\qed}
$$
\vskip-2pt

\vskip6pt
To proceed we need the following auxiliary result which will support
Remark 5.4 and Theorem 6.1$.$ It is a standard application of the
Spectral Theorem for compact operators (for similar versions see, e.g.,
\cite[Theorem 6.14.1]{NS}, \cite[Theorem 7.6]{Wei}).

\proposition {5.3}
{\it If\/ $T$ is compact, then there exist an orthonormal basis\/
$\{e_\gamma\}$ for\/ $\X$ and a family of nonnegative numbers\/
$\{\mu_\gamma\}$ such that
$$
|T|x={\sum}_\gamma\mu_\gamma\<x\,;e_\gamma\>e_\gamma
\quad\hbox{for every}\quad
x\in\X.
$$
\vskip-4pt

\proof
The operator ${|T|\in\B[\X]}$ is nonnegative (so normal) and compact$.$ (As
the class of compact operators from $\B[\X]$ is an ideal of $\B[\X]$, the
nonnegative square root $|T|$ of the nonnegative compact $|T|^2$ is compact
since ${|T|^2=T^*T}$ is compact --- see e.g., \cite[Problem 5.62]{EOT}.)
Since $\big\|\kern.5pt|T|x\big\|=\|Tx\|$ for every ${x\in\X}$ we get
${\N\big(|T|\big)=\N(T)}.$ Then by the Spectral Theorem there is a countable
orthonormal basis $\{e_k\}$ for the separable Hilbert space ${\H=\N(T)^\perp}$
consisting of eigenvectors of $|T|$ associated with positive eigenvalues
$\{\mu_k\}$ of $|T|$ such that $|T|u={\sum_k\mu_k\<u\,;e_k\>e_k}$ for every
${u\in\H}$ (see, e.g., \cite[Corollary 3.4]{ST2})$.$ Since ${\X\!=\H\oplus\N}$
with ${\N=\N(T)}$, there is an orthonormal basis
${\{e_\gamma\}\kern-1pt=\kern-1pt\{e_k\}\cup\{e_\delta\}}$ for $\X.$ Here
$\{e_\delta\}$ is an orthonormal basis for the (not necessarily separable)
Hilbert space $\N$, where ${|T|v=0}$ for ${v\in\N}$ so that the above
expansion on $\H$ describes $|T|x$ for all ${x=u\oplus v}$ in the
orthogonal direct sum ${\X\!=\H\oplus\N}$ (with ${|T|e_\delta\kern-1pt=0}$,
$\,{Te_k\kern-1pt=\mu_ke_k}$, and ${\mu_\gamma\kern-1pt=0}$ if
${\gamma\ne k)}$. \qed

\remark{5.4}
Take ${T\!\in\kern-1pt\B[\X]}$ so that ${\|Tx\|}={\big\|\kern1pt|T|x\big\|}$
for every ${x\in\X}.$ If $\{e_\gamma\}$ is any orthonormal basis for $\X$,
then
${\sum_\gamma\big\<{|T|e_\gamma\,;e_\gamma}\big\>}
\le{\sum_\gamma\big\|\kern.5pt|T|e_\gamma\big\|}
={\sum_\gamma\|Te_\gamma\|}$
and so
$$
{\sum}_\gamma\|Te_\gamma\|<\infty
\;\limply\,
{\sum}_\gamma\big\<|T|e_\gamma\,;e_\gamma\big\><\infty
\iff
T\kern-1pt\in\B_1[\X].
$$
Conversely, Proposition 5.3 ensures the existence of an orthonormal basis
$\{e_\gamma\}$ for $\X$ such that ${|T|e_\gamma=\mu_\gamma e_\gamma}.$ So
${\sum_\gamma\|Te_\gamma\|}
={{\sum}_\gamma\big\|\kern.5pt|T|e_\gamma\big\|}
={\sum_\gamma\mu_\gamma}
={\sum_\gamma\big\<|T|e_\gamma\,;e_\gamma\big\>}.$
Since ${\B_1[\X]\sse\B_\infty[\X]}$ according to Theorem 5.2(c) and
Lemma 5.1(d), we get
$$
T\kern-1pt\in\B_1[\X]
\iff
{\sum}_\gamma\big\<|T|e_\gamma\,;e_\gamma\big\><\infty
\;\limply\,
{\sum}_\gamma\|Te_\gamma\|<\infty.
$$
\vskip-6pt\noi
Therefore
$$
\hbox{$T\kern-1pt\in\B_1[\X]\iff\sum_\gamma\!\|Te_\gamma\|<\infty$
{\it for some orthonormal basis}\/ $\{e_\gamma\}$.}           \leqno{\rm(a)}
$$
Again, suppose ${\sum_\gamma\big\<|T|e_\gamma;e_\gamma\big\><\infty}$, which
means $T$ lies in $\B_1[\X].$ By Theorem 5.2(b) $\,{T\kern-1pt=A\kern1ptB}$
with ${A,B\in\B_2[\X]}$ and so ${\sum_\gamma\|Ae_\gamma\|^2<\infty}$ and
${\sum_\gamma\|Be_\gamma\|^2<\infty}$ for an arbitrary orthonormal basis
$\{e_\gamma\}$ for $\X.$ Since
${2|\<Te_\gamma,e_\gamma\>|}
={2|\<A\kern1pt Be_\gamma;e_\gamma\>|}
={2|\<Be_\gamma;A^*e_\gamma\>|}
\kern-1pt\le\kern-1pt{2\|Be_\gamma\|\kern.5pt\|A^*e_\gamma\|}
\kern-1pt\le\kern-1pt{\|Be_\gamma\|^2\!+\|A^*e_\gamma\|^2}$,
we get
$2{\sum}_\gamma|\<Te_\gamma,e_\gamma\>|
\le{\sum}_\gamma\|Be_\gamma\|^2+{\sum}_\gamma\|A^*e_\gamma\|^2
={\sum}_\gamma\|Be_\gamma\|^2+{\sum}_\gamma\|Ae_\gamma\|^2
<\infty.$
Hence
$$
\hbox{$T\kern-1pt\in\B_1[\X]\;\limply\,\sum_\gamma|\<Te_k;e_\gamma\>|<\infty$
{\it for every orthonormal basis} $\{e_\gamma\}$.}            \leqno{\rm(b)}
$$
(Actually, Claim 3(i) in the proof of Theorem 5.2 has shown by a
different proof that ${T\kern-1pt\in\B_1[\X]}$ implies
${\sum_\gamma|\<Te_k;e_\gamma\>|<\|T\|_1}).$ However, the converse of (b)
fails:
$$
\hbox{$\sum_\gamma|\<Te_k;e_\gamma\>|<\infty$
{\it for every orthonormal basis} $\{e_\gamma\}$
$\;{\kern4pt\not\kern-4pt\limply\kern4pt}\,
T\kern-1pt\in\B_1[\X]$.}                                      \leqno{\rm(c)}
$$
Indeed, take a unilateral shift ${S\in\B[\X]}$ of multiplicity one on an
infinite-dimensional separable Hilbert space $\X.$ Then $S$ shifts a countable
orthonormal basis for $\X.$ Say $Se_k=e_{k+1}$ for each integer ${k>0}$ for
some orthonormal basis $\{e_k\}$ for $\X.$ Observe that ${\<Sf_k;f_k\>=0}$ for
every orthonormal basis $\{f_k\}$ for $\X.$ In fact, take any orthonormal
basis $\{f_k\}$ for $\X$ and consider the Fourier expansion of $f_k$ in terms
of $\{e_k\}$, viz., ${f_k\!=\sum_j\<f_k;e_j\>e_j}$, and so
${Sf_k\!=\sum_j\<f_k;e_j\>Se_j}.$ Then
\begin{eqnarray*}
\<Sf_k;f_k\>
&\kern-6pt=\kern-6pt&
\Big\<{\sum}_j\<f_k;e_j\>e_{j+1};{\sum}_i\<f_k;e_i\>e_i\Big\>
={\sum}_{i,j}\<f_k;e_j\>\overline{\<f_k;e_i\>}\<e_{j+1};e_i\>            \\
&\kern-6pt=\kern-6pt&
{\sum}_j\<f_k;e_j\>\overline{\<f_k;e_{j+1}\>}
={\sum}_j\<e_{j+1};f_k\>\overline{\<e_j;f_k\>}
=\<e_{j+1};e_j\>=0,
\end{eqnarray*}
by taking the Fourier expansion of each $e_k$ in terms of $\{f_k\}.$ Thus
${\sum_k|\<Sf_k;f_k\>|}=0$ for every orthonormal basis $\{f_k\}.$ But
$S$ is an isometry, so ${S^*S=I}$ and hence ${|S|=I}$, the identity on
$\X.$ Thus ${S\not\in\B_1[\X]}$ (it is not even compact)$.$

\theorem{5.5}
{\it Let\/ $\X$ be a Hilbert space$.$ The following assertions hold true}$.$
\begin{description}
\item{$\kern-12pt$\rm(a)}
{\it $\;{(\B_1[\X],\|\cdot\|_1)}$ is a Banach space}$.$
\vskip4pt
\item{$\kern-12pt$\rm(b)}
{\it $\;\B_0[\X]$ is dense in}\/ ${(\B_1[\X],\|\cdot\|_1)}.$
\end{description}

\proof
(a)
Essentially the same argument that proves completeness of
${(\ell_1,\|\cdot\|_1)}.$ Let $\{T_n\}$ be an arbitrary $\B_1[\X]$-valued
Cauchy sequence in ${(\B_1[\X],\|\cdot\|_1)}.$ Then it is a Cauchy sequence
in the Banach space ${(\B[\X],\|\cdot\|)}$ (since ${\|\cdot\|\le\|\cdot\|_1}$),
and so
$$
\|T_n-T\|\to0
\quad\hbox{for some}\quad
T\kern-1pt\in\B[\X].
$$
Recall that the product of a pair of uniformly convergent sequences of
operators converges uniformly to the product of the limits, and also that
uniform convergence is preserved both under the adjoint and under the square
root operations (see, e.g., \cite[Problems 4.46, 5.26, 5.63]{EOT})$.$ Thus
${\|T_n-T\|=\big\|\kern1pt|T_n-T|^\half\big\|^{_{\scriptstyle2}}\to0}$ implies
$$
\big\|\kern1pt|T_n-T|^\half\big\|\to0
\quad\hbox{and}\quad
\big\|\kern1pt|T_n\big|^{_{\scriptstyle\half}}
-|T|^{_{\scriptstyle\half}}\big\|\to0.
$$
Let $\{e_\gamma\}_{\gamma\in\Gamma}\!$ be any orthonormal basis for $\X.$ Thus
${\big\|\kern1pt|T_n|^\half e_\gamma\big\|^{_{\scriptstyle2}}\!
\to\big\|\kern1pt|T|^\half e_\gamma\big\|^{_{\scriptstyle2}}}\!$
and so
$$
\|T\|_1
={\sum}_{\gamma\in\Gamma}
\big\|\kern1pt|T|^\half e_\gamma\big\|^{_{\scriptstyle2}}
\le\sup_n{\sum}_{\gamma\in\Gamma}
\big\|\kern1pt|T_n|^\half e_\gamma\big\|^{_{\scriptstyle2}}
=\sup_n\|T_n\|_1<\infty,
$$
since $\{T_n\}$ is bounded in ${(\B_1[\X],\|\cdot\|_1)}$ because it is Cauchy
in ${(\B_1[\X],\|\cdot\|_1)}.$ Then ${T\kern-1pt\in\B_1[\X]}.$
(Recall that ${\|T_n\kern-1pt-T\|\to0}$ and ${T\kern-1pt\in\B_1[\X]}$ does not
imply ${\|T_n\kern-1pt-T\|_1\kern-1pt\to0}$ --- see \cite{Kub1} for all
possible implications along this line.)
Now take an
arbitrary ${\veps>0}.$ Again, since $\{T_n\}$ is Cauchy in
${(\B_1[\X],\|\cdot\|_1)}$, there exists a finite positive integer ${n_\veps}$
such that for every ${m,n\ge n_\veps}$,
$$
{\sum}_{\gamma\in J}
\big\|\kern1pt|T_n-T_m|^\half e_\gamma\big\|^{_{\scriptstyle2}}
\le{\sum}_{\gamma\in\Gamma}
\big\|\kern1pt|T_n-T_m|^\half e_\gamma\big\|^{_{\scriptstyle2}}
=\|T_n-T_m\|_1<\veps
$$
for every finite set ${J\sse\Gamma}.$ Since
${\lim_m\!\big\|\kern1pt|T_n-T_m|^\half e_\gamma\big\|^{_{\scriptstyle2}}\!
=\big\|\kern1pt|T_n-T|^\half e_\gamma\big\|^{_{\scriptstyle2}}}\!$
for each $n$ and $e_\gamma$, it follows that
$\sum_{\gamma\in J}
\big\|\kern1pt|T_n-T|^\half e_\gamma\big\|^{_{\scriptstyle2}}\!<\veps$
for every finite set ${J\sse\Gamma}$ \hbox{and so}
$$
\|T_n-T\|_1
={\sum}_{\gamma\in\Gamma}
\big\|\kern1pt|T_n-T|^\half e_\gamma\big\|^{_{\scriptstyle2}}
=\sup_J{\sum}_{\gamma\in J}
\big\|\kern1pt|T_n-T|^\half e_\gamma\big\|^{_{\scriptstyle2}}\le\veps
$$
whenever ${n\ge n_\veps}$, where the supremum is taken over all finite sets
${J\sse\Gamma}.$ This means ${\|T_n-T\|_1\to0}.$ So every Cauchy sequence in
${(\B_1[\X],\|\cdot\|_1)}$ converges in it.

\vskip6pt\noi
(b)
Recall$:$ (i) ${\B_0[\X]\sse\B_1[\X]\sse\B_\infty[\X]}$, and (ii) Hilbert
spaces have the approximation property which means $\B_0[\X]$ is dense in
${(\B_\infty[\X],\|\cdot\|)}$, and so ${\B_0[\X]}$ is dense in
${(\B_1[\X],\|\cdot\|)}.$ To verify that ${\B_0[\X]}$ is dense in
${(\B_1[\X],\|\cdot\|_1)}$ proceed as follows$.$ Let
$\{e_\gamma\}_{\gamma\in\Gamma}\!$ be any orthonormal basis for $\X$ and take
an arbitrary ${T\kern-1pt\in\B_1[\X]}.$ So
$$
\|T\|_1={\sum}_{\gamma\in\Gamma}\big\<|T|e_\gamma\,;e_\gamma\big\>
=\sup_J{\sum}_{\gamma\in J}\big\<|T|e_\gamma\,;e_\gamma\big\><\infty,
$$
where the supremum is taken over all finite sets ${J\sse\Gamma}.$ Take an
arbitrary ${\veps>0}.$ The above expression (asserting that the family
${\{\<Te_\gamma,e_\gamma\>\}_\gamma}$ is summable) ensures the existence of
a finite set ${J_\veps\sse\Gamma}$ for which
${\sum_{\gamma\in\Gamma\\J_\veps}\big\<|T|e_\gamma\,;e_\gamma\big\><\veps}.$
Then set $\X_\veps={\span\{e_\gamma\in\X\!:\gamma\in J_\veps\}}$, a
finite-dimensional subspace of $\X$, and take $T_\veps=T|_{\X_\veps}$ in
${\B_0[\X]}.$ (Indeed, ${\R(T_\veps)=T(\X_\veps)}$ is finite-dimensional since
$\X_\veps$ is.) Thus
$$
\|T-T_\veps\|_1
={\sum}_{\gamma\in\Gamma}\big\<|T-T_\veps|e_\gamma\,;e_\gamma\big>
={\sum}_{\gamma\in\Gamma\\J_\veps}\big\<|T|e_\gamma\,;e_\gamma\big><\veps,
$$
and therefore ${\B_0[\X]}$ is dense in ${(\B_1[\X],\|\cdot\|_1)}$. \qed

\vskip9pt
Either Claim 3 in the proof of Theorem 5.2 or Remark 5.4(b) ensure that the
family ${\{\<Te_\gamma,e_\gamma\>\}_\gamma}$ is summable (i.e., the series
${\{\sum_\gamma\<Te_\gamma,e_\gamma\>\}_\gamma}$ converges in
${(\FF,|\cdot|)}\kern.5pt$) for every orthonormal basis $\{e_\gamma\}$ for the
Hilbert space $\X$, and the limit does not depend on the choice of the
orthonormal basis$.$ Therefore if ${T\kern-1pt\in\B_1[\X]}$ and $\{e_\gamma\}$
is any orthonormal basis for $\X$, then set
$$
\tr(T)={\sum}_\gamma\<Te_\gamma\kern1pt;e_\gamma\>
\quad\hbox{and hence}\quad
\|T\|_1=\tr\big(|T|\big).
$$
The number ${\tr(T)\in\FF}$ is the {\it trace}\/ of ${T\!\in\kern-1pt\B_1[\X]}$
(so the terminology {\it trace-class}\/)$.$ The norm
${\|\cdot\|_1}=\tr\big(|\cdot|\big)$ on the linear space $\B_1[\X]$ of all
trace-class operators is referred to as the {\it trace norm}\/$.$ Thus
the trace-class itself can be written as
$$
\B_1[\X]=\{T\kern-1pt\in\B[\X]\!:\;\tr\big(|T|\big)<\infty\}.
$$
By Theorem 5.2 (and Claim 3 in its proof), for every ${T\!\in\B_1[\X]}$ and
${S\kern-1pt\in\kern-1pt\B[\X]}$,
$$
|\tr(T)|\le\|T\|_1,
\qquad
\tr(T^*)=\overline{\tr(T)},
\qquad
\tr(TS)=\tr(ST),
$$
$$
\big|\tr\big(S|T|\big)\big|
=\big|\tr\big(|T|S\big)\big|
\le\|S\|\kern.5pt\|T\|_1.
$$
\vskip4pt\noi
Actually, if ${T\kern-1pt\in\B_2[\X]}$ so that ${|T|^2\in\B_1[\X]}$, then
$$
\tr\big(|T^2|\big)
=\tr(T^*T)
=\|T\|^2_2
\le\|T\|^2_1.
$$
Also, since inner products are linear in the first argument,
$$
\tr(\cdot)\!:\B_1[\X]\to\FF
\:\;\hbox{is a bounded linear functional \qquad(i.e.,}
\:\:\tr(\cdot)\in\B_1[\X]^*).
$$

\remark{5.6}
As ${\tr(\cdot)\!:\B_1[\X]\to\FF}$ is bounded and linear, and according to 
Theorem 5.2(b,c) and Lemma 5.1(e), consider the function
${\<\;\,;\;\>_2\!:\B_2[\X]\x\B_2[\X]\to\FF}\kern-1pt$ \hbox{given by}
$$
\<T;S\>_2=\tr(S^*T)
\quad\hbox{for every}\quad
S,T\kern-1pt\in\B_2[\X].
$$
Equivalently, for every ${S,T\kern-1pt\in\B_2[\X]}$ and for an arbitrary
orthonormal basis $\{e_\gamma\}$,
$$
\<T;S\>_2={\sum}_\gamma\<Te_\gamma\,;Se_\gamma\>.
$$
Since $\tr(\cdot)$ is linear, $\tr(T^*)=\overline{\tr(T)}$, and
${\<T;T\>_2}={\tr(T^*T)}={\|T\|^2_2}$, the function ${\<\cdot\,;\cdot\>_2}$ is
a Hermitian symmetric sesquilinear form which induces a quadratic form$.$ In
other words, ${\<\cdot\,;\cdot\>_2}$ is an inner product on $\B_2[\X]$ that
induces the norm ${\|\cdot\|_2}$:
\vskip6pt\noi
{\narrower
${(\B_2[\X],\|\cdot\|_2)}$ is an inner product space where the norm
${\|\cdot\|_2}$ is generated by the inner product defined by\/
$\<T;S\>_2=\tr(S^*T)$ for every\/ ${S,T\kern-1pt\in\B_2[\X]}$.
\vskip0pt}
\vskip6pt\noi
Note that the Schwartz inequality for this inner product on $\B_2[\X]$ is a
straightforward consequence of the definition of the Hilbert--Schmidt norm
${\|\cdot\|_2}$ (and, of course, of the Schwartz inequality on both Hilbert
spaces, $\X$ and $\ell^2$ over $\Gamma$):
\begin{eqnarray*}
|\<T;S\>_2|
&\kern-6pt=\kern-6pt&
|\tr(S^*T)|
=\Big|{\sum}_\gamma\<Te_\gamma\,;Se_\gamma\>\Big|
\le{\sum}_\gamma\|Te_\gamma\|\,\|Se_\gamma\|                              \\
&\kern-6pt\le\kern-6pt&
\Big({\sum}_\gamma\|Te_\gamma\|^2\Big)^{_{\scriptstyle\half}}
\Big({\sum}_\gamma\|Se_\gamma\|^2\Big)^{_{\scriptstyle\half}}\!
=\|T\|_2\,\|S\|_2
\;\;\;\hbox{for every}\;\;\;
S,T\kern-1pt\in\B_2[\X].
\end{eqnarray*}
Similarly to Theorem 5.5, it can be verified that
\begin{description}
\item{$\kern-12pt$\rm(a)}
$\;{(\B_2[\X],\|\cdot\|_2)}$ is a Hilbert space,
\vskip4pt
\item{$\kern-12pt$\rm(b)}
$\;\B_0[\X]$ is dense in ${(\B_2[\X],\|\cdot\|_2)}.$
\end{description}
In fact, a proof of (a) follows the same argument as in the proof of
Theorem 5.5(a) with $\B_1[\X]$ and${\|\cdot\|_1}$ replaced by $\B_2[\X]$
and ${\|\cdot\|_2}.$ This shows that ${(\B_2[\X],\|\cdot\|_2)}$ is complete,
thus a Banach space, and so a Hilbert space since the norm ${\|\cdot\|_2}$ is
induced by an inner product:
${\|T\|_2
={\<T;T\>_2}^{_{\scriptstyle\half}}
=\tr\big(|T|^2\big)^{_{\scriptstyle\half}}}$
for every ${T\kern-1pt\in\B_2[\X]}.$ In the same way, a proof of (b) follows
exactly as the proof of Theorem 5.5(b).

\vskip9pt
Trace and Hilbert--Schmidt classes are naturally extended in a similar fashion
to classes of operators $\B_p[\X]$ for every ${p>0}$ so that
${\B_0[\X]\subset\B_p[\X]\subset\B_q[\X]\subset\B_\infty[\X]}$ for every
${p,q}$ such that ${0<p<q<\infty}$ (see, e.g., \cite{DS2,GK,Sch2,Wei})$.$
This, however, goes beyond the scope of the present paper$.$ For further
readings on trace-class see also, for instance, \cite{RS,Con}$.$ In
particular (cf$.$ \cite[Theorem VI.26]{RS} and \cite[Theorem 19.2]{Con}),
$$
\hbox{\it $\;{(\B_1[\X],\|\cdot\|_1)}$ is a the dual of
${(\B_\infty[\X],\|\cdot\|)}$
\quad $\big($i.e., ${\B_\infty[\X]^*\cong\B_1[\X]}\kern1pt\big)$},
$$
$$
\hbox{\it $\;{(\B[\X],\|\cdot\|)}$ is a the dual of
${\big(\B_1[\X],\|\cdot\|_1)}$
\quad $($i.e\/., ${\B_1[\X]^*\cong\B[\X]}\kern1pt\big)$},
$$
which to some extent mirror the well-known classical duals ${c_0^*=\ell_1}$
and ${\ell_1^*=\ell_\infty}$ (see, e.g., \cite[Examples 1.10.3.4]{Meg} among
many others).

\section{Conclusion}

Quite often the term nuclear operator is tacitly attributed to trace class
operators without further explanation such as, for instance, ``an operator
will be called {\it nuclear}\/ if it belongs to $\B_1[\X]$''
\cite[Section III.8]{GK}$.$ This prompts our final result.

\theorem{6.1}
{\it If\/ $\X$ is a Hilbert space}\/, then
$$
\B_1[\X]=\B_N[\X]
\quad\hbox{\it and}\quad
\|\cdot\|_1=\|\cdot\|_N.
$$
\vskip-2pt

\proof
The proof is based on the he polar decomposition for Hilbert-space operators
as follows$.$ First, Proposition 5.3 is applied to show that
${\B_1[\X]\sse\B_N[\X]}$ in part (a)$.$ Then the injectiveness of $T$ when it
acts in $\N(A)^\perp$ is explored to obtain the reverse inclusion
${\B_N[\X]\sse\B_1[\X]}$ in part (b), which also yields the norm identity
${\|\cdot\|_1}{=\|\cdot\|_N}.$ Thus, to begin with, take an operator
${T\kern-1pt\in\B[\X]}$ on a Hilbert \hbox{space $\X$.}

\vskip6pt\noi
(a)
If ${T\in\B_1[\X]}$, then it is compact$.$ Thus consider the setup in the
proof of Proposition 5.3, where $\{e_\gamma\}$ is an orthonormal basis for the
Hilbert space $\X$ and $\{e_k\}$ is a countable subset of it which is an
orthonormal basis for the separable Hilbert space $\H={\N(A)^\perp\sse\X}.$
Take the direct sums $T\!={T|_\H\oplus O}$ and
$|T|={|T|\kern.5pt\big|_\H\oplus O}$ on
${\X\!=\H\oplus\N}$, where ${T|_\H\in\B[\H]}$ is injective$.$ Thus $T|_\H$
has a polar decomposition ${T|_\H=V|T|\kern.5pt\big|_\H}$ where
${V\!\in\B[\H]}$ is an isometry (see, e.g,. \cite[Corollary 5.90]{EOT}) so
that ${T=V|T|\kern.5pt\big|_\H\oplus O\in\B[\X]}={\B[\H\oplus\N]}.$ For each
$e_k$ in $\H$ set ${f_k=Ve_k}$ in $\H$ so that $\{f_k\}$ is an $\H$-valued
orthonormal sequence$.$ Then by Proposition 5.3
$$
|T|x={\sum}_\gamma\mu_\gamma\<x\,;e_\gamma\>e_\gamma
\quad\hbox{and so}\quad
Tx={\sum}_\gamma\mu_\gamma\<x\,;e_\gamma\>f_\gamma
\quad\hbox{for every}\quad
x\in\X,
$$
with ${f_\gamma=0\in\N}$ if ${\gamma\ne k}.$ Moreover, recalling that
${\mu_\gamma=0}$ if ${\gamma\ne k}$,
$$
{\sum}_\gamma\mu_\gamma
={\sum}_\gamma{\sum}_\delta\mu_\gamma\<e_\gamma\,;e_\delta\>e_\delta
={\sum}_\gamma\big\<|T|e_\gamma;e_\gamma\big\>.
$$
Therefore if ${T\kern-1pt\in\B_1[\X]}$, then
$$
{\sum}_\gamma\mu_\gamma=\|T\|_1<\infty
\quad\hbox{and}\quad
Tx={\sum}_k\mu_k\<x\,;e_k\>f_k
\quad\hbox{for every}\quad
x\in\X,
$$
where $\{e_k\}$ and $\{f_k\}$ are $\H$-valued unit sequences (thus
$\X$-valued unit sequences).

\vskip6pt\noi
(b)
Conversely, since ${\X\!=\N(T)^\perp\oplus\N(T)}$, we may regard only the
action of $T={T|_\H\oplus O}$ and ${|T|=|T|\kern.5pt\big|_\H\oplus O}$ on
${\H=\N(T)^\perp}.$ So for notational simplicity write $T$ and $|T|$ for the
injective operators $T|_\H$ and $|T|\kern.5pt\big|_\H.$ Suppose
$$
Tx={\sum}_ k\alpha_k\<x\,;z_k\>y_k
\quad\hbox{for every}\quad
x\in\H,
$$
for some $\H$-valued sequences $\{z_k\}$ and $\{y_k\}$ and for some scalar
sequence $\{\alpha_k\}$ with $\|z_k\|=\|y_k\|=1$ and
${\sum_k|\alpha_k|<\infty}.$ Thus (by polar decomposition)
$$
|T|x=V^*Tx={\sum}_ k\alpha_k\<x\,;z_k\>w_k
\quad\hbox{for every}\quad
x\in\H=\N(T)^\perp\sse\X,
$$
where ${w_k=V^*y_k}$ so that ${\|w_k\|\le1}.$ For every orthonormal basis
$\{e_j\}$ for $\H$,
$$
{\sum}_j\big\<|T|e_j\;e_j\big\>
={\sum}_j\Big\<{\sum}_k\alpha_k\<e_j;z_k\>w_k;e_j\Big\>
={\sum}_j{\sum}_k\alpha_k\<e_j;z_k\>\<w_k;e_j\>
$$
\vskip-2pt\noi
$$
\le{\sum}_k|\alpha_k|\,\Big|{\sum}_j\<e_j;z_k\>\<w_k;e_j\>\Big|
={\sum}_k|\alpha_k|\,|\<w_k;z_k\>|
\le{\sum}_k|\alpha_k|<\infty.
$$
\vskip2pt\noi
So $T$ lies in $\B_1[\H].$ Moreover, since
$Tx={\sum_ k\alpha_k\<x\,;z_k\>y_k}$ for every ${x\in\H}$,
$$
{\sum}_k\mu_k
=\|T\|_1
={\sum}_k\big\<|T|e_k;e_k\big\>
\le{\sum}_k\big\|\kern1pt|T|e_k\big\|
={\sum}_k\|Te_k\|
\le{\sum}_k|\alpha_k|,
$$
for any orthonormal basis $\{e_k\}$ for $\H.$ So
${\|T\|_1=\min\sum_k|\alpha_k|}$, the minimum taken over all scalar summable
sequences $\{\alpha_k\}$ as in the above representation of $T$$.$

\vskip6pt\noi
{From} (a) and (b) we get the following statement:
\vskip6pt\noi
{\narrower
An operator $T$ lies in $\B_1[\X]$ if and only if there are $\X$-valued
unit sequences $\{z_k\}$ and $\{y_k\}$ (i.e., ${\|z_k\|=\|y_k\|=1}$) and a
scalar summable sequence $\{\alpha_k\}$ (i.e., ${\sum_k|\alpha_k|<\infty})$
such that ${Tx=\sum_k\alpha_k\<x\,;z_k\>y_k}$ for every $x$ in $\X.$ Moreover,
$\|T\|_1={\inf\sum_k|\alpha_k]}$, where the infimum is taken over all
scalar summable sequences for which the above expression for $Tx$ holds.
\vskip0pt}
\vskip6pt\noi
Such an expression for $Tx$ is precisely a nuclear representation of $T$ as
in Section 4$.$ Then ${\B_1[\X]}={\B_N[\X]}.$ Also, with the infimum taken
over all nuclear representations of ${T\kern-1pt\in\B_N[\X]}$, we get (for
arbitrary unit sequences $\{z_k\}$ and $\{y_k\}$ and summable scalar sequence
$\{\alpha_k\}$) $\;\|T\|_N={\inf\sum_k|\alpha_k|\|z_k\|\kern1pt\|y_k\|}
={\inf\sum_k|\alpha_k|}=\|T\|_1$. \qed

\vskip6pt
Thus the notions of nuclear and trace-class coincide on Hilbert spaces$.$ For
their relationship beyond Hilbert spaces, the same first lines of Bartle's
review in Mathematical Reviews we borrowed to open the paper can be used to
close it, viz., ''Grothendieck showed that a Banach space $\X$ has the
approximation property if and only if for every nuclear operator
${T\!=\!\sum_kf_k(\cdot)y_k}$ the number ${\tr(T)\!=\!\sum_kf_k(y_k)}$ is
well-defined'' (see, e.g.,
\cite[Theorems 1.3.6, 1.3.11, 1.4.18 and Proposition 1.4.19]{DFS})$.$ This
can be regarded as a starting point for characterizing the trace property in
\hbox{Banach} spaces$.$ For further readings along this line see, for
instance, \cite{Pis1,JS,FJ}.

\section*{Acknowledgment}

The author thanks the referees for their suggestions to improve the paper.

\bibliographystyle{amsplain}

\end{document}